\documentclass[a4paper,12pt,francais]{bourbaki}
\usepackage{graphics,amscd}
\usepackage[francais]{babel}

\usepackage[matrix,arrow]{xy}

\date{Juin 2005}
\bbkannee{57\`eme ann\'ee, 2004-2005}
\bbknumero{949}
\title{ALG\`EBRES SIMPLES CENTRALES\\ 
SUR LES CORPS DE FONCTIONS DE DEUX VARIABLES }
\subtitle{[d'apr\`es A. J. de Jong]}
\author{Jean-Louis COLLIOT-TH\'EL\`ENE}
\address{Universit\'e Paris XI\\
UMR 8628 du CNRS\\
D\'epartement de Math\'ematiques\\
B\^atiment 425\\
F--91405 ORSAY C\'edex}
\email{colliot@math.u-psud.fr}
\begin{document}
\maketitle

\def\Tr  {{\rm Tr_{{\rm red}}}}
\def\ra{\rightarrow}
\def\la{\leftarrow}
\def\ind{{\rm ind}}
\def\Az{{\rm Az}}

\def\G{{\mathbb G}}
\def\A{{\mathbb A}}
\def\C{{\bf C}}
\def\R{{\bf R}}
\def\Q{{\overline Q}}
\def\Z{{\bf Z}}
\def\br{{\rm Br}}
\def\Br{{\rm Br}}
\def\pic{{\rm Pic}}
\def\Pic{{\rm Pic}}
\def\O{{\mathcal O}}
\def\Div{{\rm Div}}
\def \Ker{{\rm Ker}}
\def\P{{\mathbb P}}
\def\cl{{\rm cl}}
\def\Spec{{\rm Spec}}
\def\End{ \mathcal End  }
\def\Hom{\mathcal Hom}
\def\F{{\overline F}}

\def\V{{\overline V}}
\def\L{{\mathcal L}}

\def\et{{\rm et}}
\def\det{{\rm det}}
\def\rang{{\rm rang}}
\def\Sym{{\rm Sym}}
\def\cd{{\rm cd}}
\def\ad{{\rm ad}}
\def\ab{{\rm ab}}
\def\nr{{\rm nr}}

\def\cqfd{\hfill $\square$}
\def\Gal{{\rm Gal}}

\section{Introduction}
\`A toute alg\`ebre simple centrale $A$  de dimension finie sur
un corps $F$ sont associ\'es deux entiers,
l'indice  et l'exposant.
 L'indice de $A$ est le minimum des degr\'es $[E:F]$, o\`u $E$ 
parcourt les corps, extensions finies de $F$, pour lesquels
la $E$-alg\`ebre  $A \otimes_FE$ est isomorphe \`a une alg\`ebre de
matrices.
L'exposant de
$A$ est l'ordre de la classe de
$A$ dans le groupe de Brauer du corps $F$. L'exposant divise l'indice,
mais ne lui est pas n\'ecessairement \'egal.
Lorsque $F$ est un corps
de nombres, c'est un th\'eor\`eme des ann\'ees 1930  qu'exposant et indice
co\"{\i}ncident (\cite{[BHN]}).
  A.  J. de Jong \cite{[dJ]} montre qu'ils  co\"{\i}ncident aussi lorsque $F$ est un
corps de fonctions de deux variables sur  le corps des complexes.

Au \S 1 de cet expos\'e, on trouvera des rappels sur les alg\`ebres simples centrales
sur  un corps, le groupe de Brauer, les extensions de ces notions
au-dessus d'un sch\'ema. On passe en revue les  propri\'et\'es
connues des  corps de d\'eploiement des alg\`ebres  simples centrales et on cite
plusieurs probl\`emes ouverts dans cette direction, particuli\`erement
en ce qui concerne la  relation entre exposant et indice.

Le \S  2 d\'ecrit l'essentiel de la d\'emonstration de de Jong,
qui passe par des
d\'eformations d'alg\`ebres d'Azumaya sur une surface, rendues possibles
par  des transformations \'el\'e\-men\-taires. On d\'ecrit la d\'emonstration
{\og simplifi\'ee \fg}       \'evoqu\'ee \`a la fin de l'article de de Jong. Cette
d\'emonstration requiert un bon th\'eor\`eme de type Bertini pour
d\'emarrer, mais elle
\'evite les \'etudes fines de vari\'et\'es singuli\`eres de \cite{[dJ]}.

Le manque de temps,  de place, et de comp\'etence   nous ont fait renoncer
\`a d\'ecrire la nouvelle et tr\`es int\'eressante preuve, 
par Lieblich \cite{[Lie]},  du th\'eor\`eme de
de Jong dans le cas non ramifi\'e. Cette preuve, qui se place dans le cadre
des champs alg\'ebriques,  s'appuie
sur le th\'eor\`eme de Graber, Harris, Starr et de Jong (\cite{[GHS]},  \cite{[dJS1]})
selon lequel une vari\'et\'e (s\'eparablement) rationnellement connexe sur un corps de
fonctions d'une variable sur un corps alg\'ebriquement clos a un point rationnel.
Il est par ailleurs trop t\^ot pour parler d'une troisi\`eme d\'emonstration,
annonc\'ee par de Jong et Starr, et qui s'inscrit dans l'\'etude toute nouvelle
des vari\'et\'es 1-rationnellement  connexes.

Au \S  3, on donne une liste de propri\'et\'es des groupes alg\'ebriques
lin\'eaires sur un corps de fonctions de deux variables sur le corps
des complexes que la co\"{\i}ncidence entre exposant et indice  permet d'obtenir.

Les corps de fonctions  de deux variables sur le corps
des complexes sont de dimension cohomologique 2.
 Au \S  4, on discute divers r\'esultats r\'ecents sur des corps
qui sont de dimension cohomologique 3, \`a savoir les corps
de fonctions d'une variable sur un corps $p$-adique.

Pour de nombreuses discussions sur le travail de de Jong, je remercie
M. Ojanguren et R. Parimala. Je me suis en particulier fortement inspir\'e
d'un cours donn\'e  par cette derni\`ere \`a Lens en juin 2004. Je remercie  
 O. Gabber, P. Gille, B. Kahn, T. Szamuely et tout particuli\`erement L. Moret-Bailly
 pour leurs  critiques d'une premi\`ere version   du pr\'esent texte.

\section{La th\'eorie classique}

\subsection{ Alg\`ebres simples centrales sur un corps   (\cite{[A]},
 \cite{[D]},  \cite{[Bki]},  \cite{[Bld]}, \cite{[GSz]})}

Soient $k$ un corps et $V$ un $k$-espace vectoriel
de dimension finie $n\geq 1$. La $k$-alg\`ebre
${\rm End}(V)$ satisfait les propri\'et\'es suivantes :

\medskip

(a) Elle est de dimension finie, $n^2$, sur $k$. 

(b) Son centre est $k$.

(c) Elle ne poss\`ede pas d'id\'eal bilat\`ere non trivial.

Par le choix d'une base de $V$, ${\rm End}_k(V)$ s'identifie \`a
l'alg\`ebre
$M_n(k)$ des matrices carr\'ees de dimension $n$ sur le corps $k$.

\medskip

Une $k$-alg\`ebre simple centrale est une forme tordue
d'une telle alg\`ebre. Plus pr\'ecis\'ement :

\begin{defi}[]\label{  }Une $k$-alg\`ebre simple centrale est une
$k$-alg\`ebre de dimension finie, de centre r\'eduit \`a $k$,
sans id\'eal bilat\`ere non trivial.
\end{defi}

Pour une $k$-alg\`ebre $A$, les propri\'et\'es suivantes sont
\'equivalentes :
\medskip

(i) $A$ est une $k$-alg\`ebre simple centrale.

(ii) Il existe un corps $K$ contenant $k$ tel que la $K$-alg\`ebre
$A \otimes_kK$ est simple centrale.

(iii) Pour tout corps $K$ contenant $k$, la $K$-alg\`ebre
$A \otimes_kK$ est simple centrale.

(iv) Il existe un corps $K$ contenant $k$ et un entier $n \geq 1$ tels que
la $K$-alg\`ebre $A \otimes_kK$ est $K$-isomorphe \`a une
$K$-alg\`ebre $M_n(K)$.

(v) Il existe un corps $K$  extension  finie s\'eparable 
de $k$  et un entier $n \geq 1$  tels
 que la $K$-alg\`ebre $A \otimes_kK$ est $K$-isomorphe \`a une
$K$-alg\`ebre $M_n(K)$.

\medskip

Ceci implique en particulier que la dimension de $A$ sur $k$
est un carr\'e, soit $n^2$. L'entier $n$ est appel\'e le
{\it degr\'e} de $A$ (sur $k$).

Un corps $K$ comme en (iv) est dit {\it corps de d\'eploiement}
de l'alg\`ebre $A$. D'apr\`es (v), une cl\^oture s\'eparable
de $k$ est un corps de d\'eploiement de $A$.

D'apr\`es Wedderburn, toute $k$-alg\`ebre simple centrale $A$ est
isomorphe \`a une $k$-alg\`ebre de matrices $M_r(D)$, o\`u
$D$ est simple centrale \`a division, c'est-\`a-dire
que c'est un corps gauche de centre $k$ (on emploiera indiff\'eremment
l'une ou l'autre terminologie).
Un tel corps gauche $D$
est d\'etermin\'e, comme $k$-alg\`ebre,
 \`a isomorphisme non unique pr\`es.

 Le degr\'e de $D$ sur $k$ est
appel\'e l'{\it indice} de $A$ (sur $k$). Il est not\'e $\ind_k(A)=\ind_k(D)$.

 Soit $A$ une $k$-alg\`ebre simple centrale de degr\'e $n$.
 Pour tout  entier $m$ avec $1 \leq m \leq n$  on peut consid\'erer
  la $k$-vari\'et\'e
 alg\'ebrique dont les points sur une cl\^oture s\'eparable
 $k_s$ de $k$ sont les id\'eaux \`a droite de $A\otimes_kk_s$ 
 de $k_s$-dimension $mn$. On note cette vari\'et\'e $SB(A,m)$.
 La $k$-vari\'et\'e  $SB(A,1)$ est la vari\'et\'e de Severi-Brauer
 associ\'ee par F. Ch\^atelet \`a la $k$-alg\`ebre simple centrale $A$.
On a les faits suivants  (\cite{[Blt]}).
  La $k$-vari\'et\'e  $SB(A,m)$ est une forme de la  grassmannienne ${\rm Grass}(m,n)$.
 La $k$-vari\'et\'e $SB(A,m)$ poss\`ede un $k$-point si et seulement 
 $\ind_k(A)$ divise $m$. 
Ainsi  l'indice  $\ind_k(A)$ est le plus petit entier $m$
 tel que la $k$-vari\'et\'e $SB(A,m)$ poss\`ede un $k$-point.
 C'est aussi le p.g.c.d. des entiers $m$ satisfaisant cette propri\'et\'e.

Pour $A=M_r(D)$ comme ci-dessus, et $K$ un corps contenant $k$,
de degr\'e fini sur $k$,
les propri\'et\'es suivantes sont
\'equivalentes :  
\medskip

(i) \^Etre $k$-isomorphe \`a un sous-corps commutatif maximal de $D$.

(ii) \^Etre  de degr\'e minimal sur $k$ parmi
toutes les extensions finies de $k$  d\'eployant $A$.

\medskip

 Les sous-corps commutatifs
maximaux de $D$ sont tous de degr\'e $\ind_k(D)$. Il existe de
tels sous-corps qui sont s\'eparables sur $k$.

L'indice $\ind_k(A)$ d'une $k$-alg\`ebre simple centrale $A$
peut donc aussi se d\'efinir comme le degr\'e commun
des corps satisfaisant  l'une des propri\'et\'es ci-dessus.
C'est aussi le p.g.c.d. des degr\'es $[K:k]$
des extensions finies de corps $K/k$  d\'eployant $D$.

Une $k$-alg\`ebre simple centrale \`a division $D/k$ de degr\'e $n=\prod_i
p_i^{n_i}$
avec les $p_i$ premiers distincts s'\'ecrit comme
un produit tensoriel $D = D_1 \otimes_k \dots \otimes_k D_n$, chaque  
$D_i$
\'etant un corps gauche de degr\'e $p_i^{n_i}$.

Rappelons la notion d'alg\`ebre cyclique.
Soit $K/k$ une extension finie cyclique du corps $k$,  de degr\'e $n$, soit
$\sigma$ un g\'en\'erateur de ${\rm Gal}(K/k)$ et soit $b \in k^*$.
On munit le $k$-vectoriel $\oplus_{i=0}^{n-1}K.Y^{i}$  
d'une structure de $k$-alg\`ebre simple centrale par les relations
$Y^n=b$ et $Yx=\sigma(x)Y$ pour $x \in K$.
On note   $A=(K/k,\sigma,b)$ cette $k$-alg\`ebre.
Lorsque $k$ contient une racine $n$-i\`eme de 1, soit $\zeta_{n}$, on peut
\'ecrire $K=k(\alpha)$ avec $\alpha^n=a \in k^*$. 
Soit $\sigma$ tel que $\sigma(\alpha)=\zeta_{n}\alpha$. On note
alors $A=(a,b)_{\zeta_{n}}$. Cette $k$-alg\`ebre, consid\'er\'ee par Dickson,
 est engendr\'ee par
deux \'el\'ements $X$ et $Y$ soumis aux relations $X^n=a$, $Y^n=b$,
$YX=\zeta_nXY$.
Pour $n=2$, on retrouve la d\'efinition
des alg\`ebres de quaternions $(a,b)$.

\subsection{Groupe de Brauer d'un corps (\cite{[A]}, \cite{[Bki]},
\cite{[GSz]}, \cite{[S2]})}

Soient $A$ et $B$ deux $k$-alg\`ebres simples centrales.
Le produit tensoriel $A\otimes_kB$ est une $k$-alg\`ebre
simple centrale.
Etant donn\'ee une $k$-alg\`ebre simple centrale $A$,
on dispose de l'alg\`ebre oppos\'ee $A^{op}$, et
l'homomorphisme
$A \otimes_kA^{op} \to {\rm End}_{k-{\rm vect}}(A)$ qui \`a $a\otimes b$
associe $x \mapsto axb$ est un isomorphisme de $k$-alg\`ebres.
Consid\'erons l'ensemble des classes d'isomorphie
de $k$-alg\`ebres simples centrales. Si l'on introduit la relation
d'\'equivalence 
{\og  L'alg\`ebre $A$ est \'equivalente \`a l'alg\`ebre  $B$ s'il existe des entiers $r,s$ avec
$M_r(A) \simeq M_s(B)$ \fg},
on voit que le produit tensoriel induit
sur les classes d'\'equivalence une structure de groupe ab\'elien, dont
l'\'el\'ement neutre est la classe des alg\`ebres $M_n(k)$ ($n$
arbitraire). 
C'est le groupe de Brauer
$\Br(k)$ du corps $k$.

\medskip

Notons $\Az_n(k)$ l'ensemble des classes d'isomorphie
de $k$-alg\`ebres simples centrales de degr\'e $n$.
On dispose d'une application naturelle
$\Az_n(k) \to \Br(k)$. Cette application est une injection~:
si deux $k$-alg\`ebres simples centrales $A$ et $B$
de m\^eme degr\'e
ont m\^eme classe dans le groupe de Brauer, elles sont
isomorphes.
On a donc :
$$\Br(k) = \cup_{n=1}^{\infty} \Az_n(k),$$
la loi de groupe \'etant induite par les produits tensoriels
   $$\Az_n(k) \times \Az_m(k) \to \Az_{nm}(k).$$

\medskip

On a une seconde d\'efinition du groupe de Brauer du corps  $k$ (\cite{[S1]}).
On note $k_s$ une cl\^oture s\'eparable de $k$,
$g$ le groupe de Galois de $k_s$ sur $k$. Alors
$$ \Br(k) = H^2(g,k_s^*).$$ On passe de l'une \`a l'autre
d\'efinition en utilisant la cohomologie galoisienne
de la suite exacte de $k$-groupes lisses :
$$ 1 \to \G_{m,k} \to GL_{n,k} \to PGL_{n,k} \to 1.$$
L'ensemble point\'e de cohomologie galoisienne $H^1(k,PGL_n)=H^1(g,PGL_n(k_{s}))$
est en bijection avec $\Az_n(k)$
(\cite{[S1]}, Chap. X, \S 4 et \S 5).

\medskip

Une $k$-alg\`ebre simple centrale $A$ est d\'eploy\'ee, c'est-\`a-dire
$k$-isomorphe \`a une alg\`ebre de matrices sur $k$, 
si et seulement si sa classe $\alpha=[A] \in \Br(k)$ est nulle. Ainsi la
$K$-alg\`ebre
$A \otimes_kK$ est d\'eploy\'ee si et seulement si $\alpha_K=0 \in \Br(K)$.

\medskip

L'{\it indice}
$\ind_k(A)$ d'une $k$-alg\`ebre simple centrale $A$
ne d\'epend que de la classe $\alpha=[A]$ de $A$ dans $\Br(k)$.
Cet entier qu'on peut donc  noter $\ind_k(\alpha)$
est donc

{\rm (i)} le plus petit degr\'e d'une extension finie (s\'eparable) de corps
$K/k$ telle que l'on ait $\alpha_K=0 \in \Br(K)$;

{\rm (ii)}   le p.g.c.d.
des degr\'es des extensions finies de corps $K/k$
(s\'eparables)  telles que l'on ait
 $\alpha_K=0 \in \Br(K)$.

\medskip

On d\'efinit par ailleurs l'{\it exposant} d'une $k$-alg\`ebre simple
centrale $A$ comme l'exposant de $\alpha=[A]$ dans le groupe de
Brauer de $k$.

On utilise traditionnellement le mot  {\og exposant \fg}  (ou parfois {\og
p\'eriode \fg})
 plut\^ot que le mot
{\og ordre \fg} pour \'eviter la confusion avec les ordres
(maximaux ou autres) lorsque le corps $k$ est le corps
des fractions d'un anneau.

\bigskip

\begin{prop}[R. Brauer]\label{}

{\rm (i)} L'exposant divise l'indice.

{\rm (ii)} Les nombres premiers qui divisent l'indice divisent l'exposant.
\end{prop}
 
 \noindent{\sc Preuve} ---
Le premier \'enonc\'e r\'esulte de l'existence, pour une extension
finie de corps $K/k$, d'un homomorphisme  de corestriction $\Br(K) \to
\Br(k)$
pour lequel la composition avec la restriction $\Br(k) \to \Br(K)$
est la multiplication par le degr\'e $[K:k]$.

Montrons le second \'enonc\'e. Soit $l$ premier ne divisant pas
l'exposant
 de $A$. Soit $K/k$ une extension finie galoisienne d\'eployant
$A$. Soit $F \subset K$ le corps fixe d'un $l$-sous-groupe de Sylow
de $\Gal(K/k)$. L'exposant de $[A\otimes_kF] \in \Br(F)$ divise celui de
$[A] \in \Br(k)$, il est donc premier \`a $l$. Par ailleurs
la restriction de $A\otimes_kF$ \`a $K$ est triviale, l'argument
de corestriction montre que la classe $[A\otimes_kF] \in \Br(F)$
est annul\'ee par $[K:F]$ qui est une puissance de $l$. Ainsi
$[A\otimes_kF]=0 \in \Br(F)$, la $k$-alg\`ebre $A$ est d\'eploy\'ee
par l'extension $F/k$ qui est de degr\'e premier \`a $l$. Ainsi
$l$ ne divise pas l'indice. \cqfd

\subsection{Alg\`ebres d'Azumaya, cohomologie \'etale, 
ramification (\cite{[Gr]}, \cite{[Mi]})}

Nous nous contenterons ici de rappeler quelques r\'esultats. 
Pour les d\'emonstrations, on renvoie aux trois expos\'es
de Grothendieck \cite{[Gr]} et au chapitre IV du livre de Milne
\cite{[Mi]}.

 Soit $X$ un sch\'ema. Une alg\`ebre d'Azumaya sur $X$ de degr\'e
$n$ est un faisceau de $\O_X$-alg\`ebres 
localement libres de rang fini 
qui localement pour la topologie
\'etale sur $X$ est isomorphe \`a $M_n(\O_X)$.

Soit $A$ une telle alg\`ebre. \`A tout entier $m$ avec $1 \leq m \leq n$ on associe
un $X$-sch\'ema $SB(A,m)$. C'est le sch\'ema des
$\O_X$-id\'eaux \`a droite de $A$ qui sont localement libres de
rang $mn$ sur $X$  et qui sont localement facteurs directs dans $A$.
C'est un $X$-sch\'ema projectif, lisse, \`a fibres connexes.
Pour $m=1$, c'est le sch\'ema de Severi-Brauer (\cite{[Gr]})  associ\'e \`a $A$. 

 On notera $\Az_n(X)$ l'ensemble des classes d'isomorphie d'alg\`ebres
 d'Azumaya sur $X$ de degr\'e $n$.
Le produit tensoriel de telles $\O_X$-alg\`ebres induit un
produit
$$\Az_n(X)\times \Az_m(X) \to \Az_{n+m}(X).$$
Si l'on consid\`ere l'application induite sur les
classes d'isomorphie, et que de plus on consid\`ere comme triviales  les
alg\`ebres de la forme
$\End(V)$ pour
$V$ un fibr\'e vectoriel sur $X$, on obtient un groupe ab\'elien
 $\br_{\Az}(X)$ (ce groupe est not\'e $\br(X)$ dans \cite{[Gr]}). 
 Pour $X$ quasi-compact, 
 ce groupe est de torsion.

On dispose par ailleurs du groupe de Brauer d\'efini par Grothendieck.
C'est le deuxi\`eme groupe de cohomologie \'etale
$\br(X)= H^2(X,\G_m)$\ (ce groupe est not\'e $\br'(X)$ dans \cite{[Gr]}). 

Si le sch\'ema $X$ est r\'egulier, ou plus g\'en\'eralement
si tout ouvert \'etale de $X$ a ses anneaux locaux factoriels,
le groupe $\br(X)$ est de torsion : ceci r\' esulte du th\'eor\`eme
\ref{injec} ci-dessous.

Pour tout sch\'ema $X$, il y a un plongement naturel
$$\br_{\Az}(X) \hookrightarrow \br(X).$$
Un th\'eor\`eme de Gabber, dont de Jong [dJ2] a donn\'e une
d\'emonstration,
assure que si $X$ est quasi-compact et quasi-s\'epar\'e,
et poss\`ede un fibr\'e inversible ample, alors
le plongement ci-dessus induit un isomorphisme entre
$\br_{\Az}(X)$ et le sous-groupe de torsion de $\br(X)$.

\bigskip

Soit $R$ un anneau de valuation discr\`ete de corps des fractions $K$,
de corps r\'esiduel $\kappa$. Soit $n$ un entier inversible sur $R$.
On dispose alors d'une suite exacte naturelle
$$ 0 \to {}_n\br(R) \to {}_n\br(K) \to H^1(\kappa,\Z/n) \to 0.$$
(Pour un groupe ab\'elien $M$  et un entier $n>0$,
on note ${}_nM=\{x \in M, nx=0.\}$.)

On note $\partial_R $
l'homomorphisme ${}_n\br(K) \to H^1(\kappa,\Z/n)$, qu'on appelle
l'application r\'esidu. Un \'el\'ement de
$\alpha \in {}_n\br(K)$ est dit {\it ramifi\'e}  (par rapport  \`a
l'anneau de valuation discr\`ete $R$)
si $\partial_R(\alpha)\neq 0$.

\medskip

Soit $R \subset S$ une inclusion locale d'anneaux de valuation discr\`ete 
induisant une extension finie $K \subset L$ des corps des fractions
et une  extension finie   $\kappa_R \subset \kappa_S$ 
des corps r\'esiduels.
Soit $e=e_{S/R}$ l'indice de ramification de $S$ sur $R$. Soit $n>0$ premier
aux caract\'eristiques r\'esiduelles. On a alors le diagramme
commutatif :

\[\xymatrix{
        {}_n\Br(K) \ar[r]^{\partial_R} \ar[d]_{{\rm Res}_{K,L}} &
H^1(\kappa_R,{\bf Z}/n) 
\ar[d]^{{\times e_{S/R}.{\rm
Res}_{\kappa_R,\kappa_S}}} \\
        {}_n\Br(L) \ar[r]^{{\partial_S}} & H^1(\kappa_S,{\bf Z}/n).   }
\]

\medskip

Soit $X$ un sch\'ema noeth\' erien int\`egre.
Soit $K$ le corps des fonctions
rationnelles de $X$.
 On note  $X^{(i)}$ l'ensemble des
points de codimension $i$ de
$X$.  On note $\kappa_x$
le corps r\'esiduel d'un point $x \in X$.
On note $\mu_l$ le faisceau \'etale des racines $l$-i\`emes de 1,
et $\mu_l^{\otimes 2}=\mu_l \otimes \mu_l$.

\begin{prop}[]\label{1.3.1} Soit $X$ un sch\'ema noeth\'erien int\`egre, excellent, de dimension {\rm 2}.  Soit
$l$ un nombre premier inversible sur $X$.

{\rm (i)} Il existe un complexe naturel de groupes de cohomologie \'etale
$$0 \to H^2(X,\mu_l^{\otimes 2})  \to H^2(K,\mu_l^{\otimes 2}) \to
\oplus_{x\in X^{(1)}} H^1(\kappa_x,\mu_l)
\to  \oplus_{x\in X^{(2)}} H^0(\kappa_x, \Z/l) \to  0.$$

{\rm (ii)} Si $X$ est r\'egulier,
ce complexe est exact en le terme $H^2(K,\mu_l^{\otimes 2})$.

{\rm (iii)} Si $X$ est un sch\'ema local r\'egulier,
le complexe est aussi exact en $H^2(X,\mu_l^{\otimes 2})$.
\end{prop}

\noindent{\sc Quelques r\'ef\'erences} --- 
Pour une description du complexe, et en particulier
  des fl\`eches dans le complexe,
je renvoie \`a l'article de Kato \cite{[Kt]} (qui d\'ecrit
de tels complexes dans un cadre bien plus large que
celui utile pour le pr\'esent expos\'e).
Pour \'etablir les autres \'enonc\'es on peut
adjoindre les racines $l$-i\`emes de 1,
on est alors ramen\'e \`a deux th\'eor\`emes sur
le groupe de Brauer, dus \`a Auslander-Goldman et
Grothendieck. L'\'enonc\'e (iii) est une cons\'equence
du th\'eor\`eme \ref{purdim2} ci-apr\`es car pour $X$ local, 
soit $X=\Spec(R)$, on a $H^2(R,\mu_l)={}_lH^2(R,\G_m)$. 
L'\'enonc\'e (ii) r\'esulte du th\'eor\`eme \ref{injec} ci-apr\`es.
 (Dans le cas local r\'egulier, on conjecture que tout le
complexe est exact, ce n'est pas utile pour notre propos.)

\begin{theo}[\cite{[AG]}, \cite{[Gr]}]\label{injec}
Soit $X$ un sch\'ema  noeth\'erien int\`egre,  et soit $K$ son corps des 
fonctions. Si $X$ est  g\'eom\'etriquement localement factoriel, par
exemple si $X$ est r\'egulier, l'homomorphisme $\Br(X) \to \Br(K)$ est
injectif. \cqfd
\end{theo}

\begin{theo}[]\label{purdim2}
Soit $X$ un sch\'ema  r\'egulier int\`egre de dimension au plus $2$,
et soit $K$ son corps des fonctions.  Tout \'el\'ement de
$\Az_n(K)$ dont l'image dans $\Br(K)$ est dans
$\Br(X) \subset \Br(K)$ est dans l'image de l'application
de restriction $\Az_n(X) \to \Az_n(K)$.

Soit $n>0$  un entier inversible sur $X$. 
Si un \'el\'ement de $\Az_n(K)$ a une image dans ${}_n\Br(K)$ dont
les r\'esidus en tous les points de codimension 1 de $X$ sont nuls,
alors il provient d'un \'el\'ement de $\Az_n(X)$. 
\end{theo}
Ce th\'eor\`eme combine un \'enonc\'e sur les anneaux de
valuation discr\`ete (\cite{[R]}) et le th\'eor\`eme que tout module r\'eflexif sur un
anneau local r\'egulier de dimension $2$ est libre.
On renvoie le lecteur \`a \cite{[AG]}, \cite{[Gr]} (II, Thm. 2.1),
\cite{[CTS]} Cor. 6.14),
\cite{[OP]} (\S 1). \cqfd

\medskip

Soient $X$ et $K$ comme dans la proposition \ref{1.3.1}. \`A tout \'el\'ement $\alpha \in
H^2(K,\mu_l^{\otimes 2})$ on associe son lieu de ramification : c'est
l'adh\'erence dans $X$ de l'ensemble des points $x$ de codimension 1 de $X$
en lesquels $\alpha$ a un r\'esidu non nul dans $H^1(\kappa_x,\mu_l)$.

\`A tout couple d'\'el\'ements $u,v \in K^*$ on associe le symbole $(u,v) \in H^2(K,\mu_l^{\otimes 2})$ :
c'est le cup-produit des classes de $u$ et $v$ dans $K^*/K^{*l}=H^1(K,\mu_l)$.
Lorsque les racines $l$-i\`emes de 1 sont dans $K$, apr\`es avoir fait un choix d'une racine
primitive, on trouve ainsi les classes des alg\`ebres cycliques dans $H^2(K,\mu_l)={}_l\Br(K)$.

\begin{prop}[Saltman \cite{[Sa1]}, \cite{[CTOP]}]\label{1.3.2}{\it Soit $R$ un anneau local
r\'egulier excellent de dimension  $2$, d'id\'eal maximal $m$, de   corps
des fractions $K$. Soit $\alpha \in H^2(K,\mu_l^{\otimes 2})$
un
\'el\'ement dont le lieu de ramification sur $X=\Spec(R)$ est un diviseur
strictement \`a croisements normaux.

{\rm (i)} Si ce lieu est vide, alors  $\alpha$ est dans l'image de
$H^2(R,\mu_l^{\otimes 2})$.

{\rm (ii)} Si ce lieu est de la forme $s=0$, o\`u $s \in m$ est un
param\`etre r\'egulier de $R$ (i.e. $R/s$ r\'egulier), alors il existe $u
\in R^*$ tel que $\alpha - (u,s)$ soit dans l'image de
$H^2(R,\mu_l^{\otimes 2})$.

{\rm (iii)} Si ce lieu est de la forme $st=0$, o\`u $s$ et $t$ engendrent l'id\'eal $m$, alors
 il existe $u,v \in R^*$ et $r \in \Z/l$ tels que $\alpha-
 (u,s) - (v,t) -  r(s,t)
 \in H^2(K,\mu_l^{\otimes 2})$
soit dans l'image de $H^2(R,\mu_l^{\otimes 2})$.}
\end{prop}

\noindent{\sc Preuve} --- Le (i) r\'esulte de la proposition
\ref{1.3.1} (ii).

 Notons $\kappa_s$ le corps des fractions de $R/s$.
Le seul r\'esidu non trivial de $\alpha$ est en $s=0$,
c'est une classe  $\partial_s(\alpha)=\xi_s \in
\kappa_s^*/(\kappa_s^*)^l$
dont l'image dans $\Z/l$ par l'application
valuation est z\'ero (\ref{1.3.1} (i)). On a la suite exacte \'evidente
$$1 \to  (R/s)^*/(R/s)^{*l} \to \kappa_s^*/(\kappa_s^*)^l \to \Z/l \to
0.$$ Il existe donc un \'el\'ement $u_s
\in (R/s)^*$ qui a pour image $\xi_s$ via l'application de r\'eduction
$(R/s)^*/(R/s)^{*l} \to \kappa_s^*/\kappa_s^{*l}$.
L'application de r\'eduction $R^* \to (R/s)^*$ est surjective
($R$ est un anneau local),
on peut donc relever $u_s$ en $u \in R^*$.
Les formules usuelles pour le r\'esidu d'un cup-produit
montrent  que le r\'esidu  de $(u,s)$ en $s=0$ est $\xi_s$,
et il est z\'ero partout ailleurs, comme celui de $\alpha$.
La proposition \ref{1.3.1} (ii)
donne alors l'\'enonc\'e (ii).

Montrons (iii). Notons $\kappa_s$, resp. $\kappa_t$,
 le corps des fractions de $R/s$, resp. $R/t$, et $\kappa$ le
 corps r\'esiduel en l'id\'eal maximal $m$ de $R$.
 Notons $\partial_s : H^1(\kappa_s,\mu_l) \to H^0(\kappa,\Z/l)$
 et de m\^eme $\partial_t : H^1(\kappa_t,\mu_l) \to H^0(\kappa,\Z/l)$
 les homomorphismes apparaissant dans la proposition \ref{1.3.1}.
Soient $\xi_s$ et $\xi_t$ les r\'esidus de
$\alpha$ en $s=0$ et $t=0$. Comme la proposition
\ref{1.3.1} donne un complexe de groupes ab\'eliens,
  on a
$\partial_s(\xi_s)+\partial_s(\xi_t)=0 \in \Z/l$.
Soit $r \in \Z$ avec
$\partial_s(\xi_s)=r \in \Z/l$ et
$\partial_t(\xi_t)=-r \in \Z/l$.
Proc\'edant comme ci-dessus, on trouve $u,v \in R^*$
tels que l'image de $ut^r \in R[1/t]^*$ dans $\kappa_s^*/\kappa_s^{*l}$
soit $\xi_s$ et que l'image de $vs^{-r} \in R[1/s]^*$ dans
$\kappa_t^*/\kappa_t^{*l}$ soit $\xi_t$.
On v\'erifie   que tous les r\'esidus de
$$\alpha+(s,u)+(t,v)+r(s,t) \in H^2(K,\mu_l^{\otimes 2}) $$
aux points de codimension 1 de $X$ sont nuls,
la proposition \ref{1.3.1} (ii)  donne alors
que cet \'el\'ement de $H^2(K,\mu_l^{\otimes 2})$
est  l'image d'un (unique)  \'el\'ement  de $H^2(R,\mu_l^{\otimes 2})$.
\cqfd

\subsection{ Questions sur les corps de d\'eploiement}

Sauf mention explicite du contraire, on suppose ici que
les corps consid\'er\'es sont de carac\-t\'eris\-tique z\'ero,
ou du moins que l'indice des alg\`ebres consid\'er\'ees
est premier \`a la caract\'eristique du corps de base.
Etant donn\'es un corps $k$ et une $k$-alg\`ebre simple
centrale $A$, que peut-on dire sur les extensions finies
de corps $K/k$ qui d\'eploient $A$ ?

\bigskip

Il existe des $k$-alg\`ebres simples centrales qui ne sont pas d\'eploy\'ees
par une extension cyclique du corps de base
(Tignol-Amitsur, Tignol-Wadsworth  \cite{[TW]} Exemples 3.6 et  Thm. 4.7 (v)).
Un exemple  simple est le produit tensoriel
$(x,y)\otimes_k(z,t)$ des alg\`ebres de quaternions $(x,y)$ et $(z,t)$
sur le corps $k={\bf C}(x,y,z,t)$ et m\^eme d\'ej\`a 
sur le corps de s\'eries formelles it\'er\'ees $k=\C((x))((y))((z))((t))$.

Toute $k$-alg\`ebre simple centrale $A$ d'exposant $n$ sur un corps $k$
contenant
une racine primitive $n$-i\`eme de 1, soit $\zeta_n$, est semblable \`a un
produit d'alg\`ebres cycliques $(a,b)_{\zeta_n}$, 
elle est en particulier
d\'eploy\'ee par une extension multicyclique de $k$ 
(cons\'equence  imm\'ediate du th\'eor\`eme de Merkur'ev et Suslin \cite{[MS]}). Il en r\'esulte que
toute $k$-alg\`ebre simple centrale $A$ est
d\'eploy\'ee par une extension r\'esoluble de $k$.

\bigskip

\' Etant donn\'es un corps $K$ et une $k$-alg\`ebre simple centrale
{\it \`a division} $D$, on dit qu'elle est un produit crois\'e
si elle admet un sous-corps commutatif maximal $K \subset D$ galoisien
sur $k$. On dit qu'elle est cyclique si elle admet un sous-corps
commutatif maximal $K \subset D$ cyclique sur $k$.

\bigskip
Amitsur (1972) montra qu'il existe des alg\`ebres \`a division qui ne sont
pas des produits crois\'es : elles ne poss\`edent pas de sous-corps
commutatif maximal galoisien sur le corps de base.
Les indices de ses exemples ainsi que de ceux qui ont \'et\'e construits
par la suite
sont de la forme $2^r \prod_i p_i^{n_i}$ avec les $p_i \geq 3$ premiers,
$r \geq 3$, et $n_i \geq 2$ pour tout $i$.
Lorsque l'indice est premier,  la question suivante est ouverte :

\medskip

{\bf Probl\`eme} {\it Sur tout corps $k$, toute alg\`ebre \`a division 
d'indice premier
$l$ est-elle cyclique ?}

\medskip

C'est le cas de fa\c con triviale  pour $l=2$, c'est un r\'esultat
de Wedderburn pour $l=3$, la question est ouverte d\'ej\`a pour
$l=5$.

\bigskip

D\`es les ann\'ees 1930, des exemples furent donn\'es (par
Albert, Brauer, K\"othe, Nakayama, voir les r\'ef\'erences dans \cite{[CTG]})
pour montrer qu'indice et exposant d'une alg\`ebre
 ne co\"{\i}ncident pas n\'ecessairement.

Un \'enonc\'e g\'en\'eral permettant de fabriquer de tels exemples
est le suivant :

 \begin{prop}[Tignol \cite{[T]}] \label{tignol}
Soient $K/k$ une extension cyclique de corps,
$\sigma$ un g\'en\'erateur du groupe de Galois de $K$ sur $k$
et $t$ une variable. Soit $A$ une $k$-alg\`ebre simple centrale.
On a  les formules
$$ \ind_{k(t)}(A_{k(t)}\otimes_{k(t)}(K(t)/k(t),\sigma,t))=\ind(A_K).[K:k]$$
et
$$
\ind_{k((t))}(A_{k((t))}\otimes_{k((t))}(K((t))/k((t)),\sigma,t))=\ind(A_K).[K:k].$$
En particulier, si $A$ est un corps gauche, les conditions suivantes sont
\'equivalentes :

{\rm (i)} La $K$-alg\`ebre $A_K$ est un corps gauche.

{\rm (ii)} La $k(t)$-alg\`ebre
$A_{k(t)}\otimes_{k(t)}(K(t)/k(t),\sigma,t)$ est un corps gauche.

{\rm (iii)} La $k((t))$-alg\`ebre
$A_{k((t))}\otimes_{k((t))}(K((t))/k((t)),\sigma,t)$ est un corps
gauche.\cqfd
\end{prop}

Donnons deux exemples. 
Sur le corps $K=\C(x_{1},\dots,x_{d})$ des fonctions rationnelles
en $d$ variables sur le corps des complexes, le produit
tensoriel d'alg\`ebres de quaternions $A=(x_{1}+2, x_{2}) \otimes_{K} \dots
\otimes_{K} (x_{1}+d,x_{d})$ est un corps gauche, d'exposant $2$ et d'indice $2^{d-1}$.
Il en est ainsi d\'ej\`a sur le corps $K=\C(x_{1})((x_2)) \dots ((x_{d}))$.
Par ailleurs soit  ${\bf F}$ un corps fini  de caract\'eristique diff\'erente de 2,
et soit $a \in {\bf F}$ non carr\'e. 
Sur le corps
de fonctions rationnelles $K={\bf  F}(x,y)$, l'alg\`ebre
$(x,a)\otimes_{K}(x+1,y)$ est un corps gauche, d'exposant 2 et d'indice 4.
Ces exemples font intervenir la ramification des alg\`ebres consid\'er\'ees
en des valuations convenables du corps de base.

Dans \cite{[Kr]} et \cite{[CTG]}
on trouvera diverses constructions d'alg\`ebres d'Azumaya $A$ (donc
sans ramification)
sur des vari\'et\'es projectives, lisses, connexes $X$ sur les complexes
(de dimension au moins 3) telles que l'exposant divise proprement
l'indice de la ${\bf C}(X)$-alg\`ebre obtenue par \'evaluation
de $A$ au point g\'en\'erique de $X$. Pour tout $d\geq 2$, et tout nombre premier $l$,
 O.~Gabber donne des exemples de vari\'et\'es projectives, lisses, connexes $X$
 de dimension $d$, en fait des produits de courbes,  qui poss\`edent  une alg\`ebre d'Azumaya
 d'exposant $l$ et d'indice (g\'en\'erique) $l^{d-1}$.
\medskip

Les exemples obtenus de cette fa\c con sont d\'efinis sur des
corps dont la dimension coho\-mo\-logique \cite{[S2]} est au moins 3.
Par d'autres m\'ethodes, pour tout entier $n$,  Merkur'ev \cite{[M2]} 
construit un (tr\`es gros) corps $K_n$ de dimension cohomologique 2
sur lequel il existe un produit tensoriel de $n$
alg\`ebres de quaternions  qui est un corps gauche,  et donc
d\'efinit un \'el\'ement de $\Br(K_n)$ d'exposant 2 et d'indice $2^n$.

\medskip

Il existe n\'eanmoins plusieurs classes de  corps de dimension
cohomologique 2 ([S2]) (ou de dimension cohomologique virtuelle 2)
pour lesquels on sait que l'indice co\"{\i}ncide
avec l'exposant. Pour chacune des classes suivantes, on sait
en outre que toute alg\`ebre \`a division est cyclique.
\begin{itemize} 

\item[$\bullet$] Les corps globaux et les
corps locaux (Brauer-Hasse-Noether, Albert) (\cite{[BHN]}, \cite{[D]}, \cite{[Ro]}).

\item[$\bullet$] Les corps de la forme $k=F((t))$ avec $F$ corps de
caract\'eristique z\'ero et de dimension cohomologique 1.

\item[$\bullet$] Les
corps de fractions d'anneaux locaux normaux hens\'eliens excellents de dimension
2,
\`a corps r\'esiduel alg\'ebriquement clos de caract\'eris\-tique
z\'ero (Artin \cite{[Ar3]}, 
Ford-Saltman \cite{[FS]}, \cite{[CTOP]}). 
\end{itemize}

\medskip

Soit $r \geq 0$ un entier. Un corps $k$ satisfait la propri\'et\'e $C_r$ si
 toute forme homog\`ene \`a 
coefficients dans $k$, de degr\'e $d$,
en $n>d^r$ variables, a un z\'ero non trivial dans $K$. Les corps de fonctions de $r$ variables sur un  corps alg\' ebriquement clos
 sont des
corps $C_r$ (S.~Lang). C'est une cons\'equence
du th\'eor\`eme de Merkur'ev et Suslin \cite{[MS]} que tout corps
$C_2$ est de dimension cohomologique au plus 2.

\medskip

{\bf Probl\`eme}  (M. Artin  \cite{[Ar2]}, Appendix) {\it Soit $k$ un corps $C_2$.
Pour toute alg\`ebre simple centrale sur $k$,
l'indice est-il \'egal \`a l'exposant ?}

\medskip

Lorsque l'on se limite aux alg\`ebres d'indice 2-primaire 
ou 3-primaire, la r\'eponse est affirmative\footnote{Quitte \`a remplacer
l'hypoth\`ese $C_{2}$ par sa variante $C'_{2}$, cf.   \cite{[CTG]}, Prop.7.)}. 
 Cette cons\'equence du th\'eor\`eme de
Merkur'ev et Suslin \cite{[MS]} a \'et\'e remarqu\'ee par
plusieurs auteurs (\cite{[Ar2]}, \cite{[MS]}). Le cas particulier des corps de fonctions de deux variables
sur les complexes avait \'et\'e \'etabli plus t\^ot, par M.~Artin et J.~Tate 
 (\cite{[Ar2]}, Appendix),
qui utilisaient un r\'esultat de S.~Bloch (1974).

Pour les corps de fonctions de deux variables sur les complexes,
et les alg\`ebres d'indice arbitraire sur de tels corps, la r\'eponse 
affirmative
\`a cette question est le th\'eor\`eme de {de~Jong}  \cite{[dJ]} qui fait 
l'objet de cet expos\'e. Pour $r \leq 2$, on a donc une r\'eponse
affirmative \`a la question suivante, o\`u la borne sugg\' er\'ee
est, d'apr\`es les exemples mentionn\' es ci-dessus, la meilleure
possible.

\medskip

{\bf Probl\`eme}\footnote{Soit $X$ une vari\'et\'e de dimension $r$ sur le corps  $\R$ des r\'eels,
g\'eom\'etriquement int\`egre et sans point  r\'eel. Le corps des fonctions
$K=\R(X)$ est de dimension cohomologique $r$ mais on ne sait pas s'il est $C_r$.
On peut pour un tel corps poser le m\^eme probl\`eme que ci-dessus. 
Pour les surfaces, i.e. pour $r=2$,
le probl\`eme se ram\`ene \`a  l'une quelconque des questions suivantes (Pfister, 1982) :

(a) Toute forme quadratique en au moins 5 variables sur un tel corps admet-elle un z\'ero
non trivial   ?

(b) Le produit tensoriel de deux alg\`ebres de quaternions sur un tel corps
est-il semblable \`a une alg\`ebre de quaternions ?
}
 {\it Pour toute alg\`ebre simple centrale sur un corps de
fonctions de
$r$ variables sur le corps des complexes, l'indice divise-t-il l'exposant \`a la
puissance 
$r-1$ ?}

\medskip
Les corps de fonctions de $r$ variables sur le corps des complexes
sont des corps $C_{r}$. On peut se poser la question ci-dessus pour tout corps $C_{r}$.
Les corps de fonctions de $r$ variables sur un corps fini sont des corps $C_{r+1}$.

\medskip

 {\bf Probl\`eme}  {\it  Pour toute alg\`ebre simple centrale sur  un
corps de fonctions de $r$ variables sur un corps fini, l'indice
divise-t-il l'exposant \`a la puissance $r$ ?}

Lorsque $r=1$, c'est un r\'esultat classique (Hasse).

Lorsque $r=2$, c'est un r\'esultat r\'ecent de Lieblich \cite{[Lie]}
que  pour les alg\`ebres non ramifi\'ees, d'indice premier \`a la caract\'eristique,
 l'indice est \'egal \`a l'exposant\footnote{
Pour obtenir ce
r\'esultat, outre des techniques de champs alg\'ebriques il utilise
les propri\'et\'es des espaces de modules de fibr\'es vectoriels
sur les surfaces projectives et lisses  \cite{[HL]}.}.
Pour une alg\`ebre
ramifi\'ee, la technique \'evoqu\'ee au d\'ebut du \S 4 ci-dessous
(m\'ethode de Saltman) permet alors de montrer que
l'indice divise l'exposant au cube. Comme on
a vu ci-dessus, dans ce cas le mieux que l'on puisse  esp\'erer en g\'en\'eral est
que l'indice divise l'exposant  au carr\'e.

\section{Le th\'eor\`eme  de de Jong}

\subsection{Deux invariants cohomologiques  des
alg\`ebres d'Azumaya de degr\'e $n$}

Soit  $X$ un sch\'ema et $n>0$ un entier. L'ensemble $\Az_n(X)$
des classes d'isomorphie d'alg\`ebres d'Azumaya
de degr\'e $n$ sur $X$ s'identifie \`a l'ensemble de cohomologie de
\v{C}ech \'etale   $H^1(X,PGL_n)$.

On a le diagramme commutatif de suites exactes de $X$-sch\'emas en groupes
 suivant :

\[\xymatrix{
& 1   \ar[d]  & 1 \ar[d]  && \\
       1 \ar[r] & \mu_n \ar[r]  \ar[d]_{} & SL_n \ar[r]  \ar[d]  &
                  PGL_n \ar[r] \ar[d]^{=}  & 1 \\
       1 \ar[r] & \G_m \ar[r]\ar[d]_{x \mapsto x^n}  & GL_n \ar[r]\ar[d]^{\det} & PGL_n  \ar[r] & 1 \\
 & \G_m \ar[d]\ar[r]^{=} &\G_m \ar[d]  &   & 
\\
& 1      & 1    & .&
 }
\]

\medskip

Supposons $n$ inversible sur $X$. Les suites exactes
ci-dessus d\'efinissent alors des suites exactes de
faisceaux pour la topologie \'etale sur $X$.
De la suite horizontale m\'ediane on tire la suite exacte d'ensembles point\'es de
cohomologie  de \v{C}ech \'etale (\cite{[Mi]}, Chapitre IV, Thm. 2.5) :
$$H^1(X,\G_m) \to H^1(X,GL_n) \to H^1(X,PGL_n) \to
{}_nH^2(X,\G_m). \hskip1cm (3.1.1)$$

L'application compos\'ee $\Az_n(X)=H^1(X,PGL_n) \to
H^2(X,\G_m)=\Br(X) $ associe \`a une alg\`ebre d'Azumaya $A$ sa
classe $[A]$ dans le groupe de Brauer cohomologique $\Br(X)$,
classe qui est annul\'ee par $n$.

La suite exacte sup\'erieure donne naissance \`a une application
$$\cl : \Az_n(X)=H^1(X,PGL_n) \to H^2(X,\mu_n).$$

La suite exacte verticale de gauche (suite de Kummer)
donne naissance \`a la suite exacte bien connue
$$ 0 \to \pic(X)/n \to H^2(X,\mu_n) \to {}_n\Br(X) \to 0.$$

Pour $A \in \Az_n(X)$, la  fl\`eche de droite envoie $\cl(A)\in H^2(X,\mu_n) $ 
sur $[A] \in {}_n\Br(X)$.

  \begin{lemm}[]\label{2.1.1}
 {\it Soit $X$ un sch\'ema.

{\rm (a)} Soit $A \in \Az_n(X)$. Si l'on a $[A]=0 \in \Br(X)$, alors il existe un
fibr\'e
vectoriel $V$ sur $X$ de rang $n$ tel que $A=\End_X(V)$.

{\rm  (b)} Soient  $V_1$ et $V_2$ deux fibr\'es vectoriels  de rang $n$ sur $X$;
si l'on a un isomorphisme d'alg\`ebres $\End(V_1) \simeq \End(V_2)$,
alors il existe un fibr\'e inversible
$L$ sur $X$ et un isomorphisme  $V_1 \simeq V_2 \otimes L$.

{\rm  (c)} Soit $V$ un fibr\'e vectoriel sur $X$ de rang $n$,
$A={\End}(V) \in \Az_n(X)$. Supposons $n$ inversible sur $X$.
L'image de la classe de
$\det(V)=\Lambda^nV$ dans $\Pic(X)$ par la fl\`eche de Kummer
$ \Pic(X)/n  \to H^2(X,\mu_n)$ co\"{\i}ncide avec l'oppos\'e de $\cl(A) \in
H^2(X,\mu_n)$. \cqfd}
\end{lemm}

 \noindent{\sc Preuve} --- Les  points (a) et (b)  sont des cons\'equences de
 la suite exacte (3.1.1).  Le point~Ê(c) s'\'etablit en consid\'erant les suites exactes de cohomologie 
de \v{C}ech \'etales d\'eduites du diagramme commutatif de suites exactes
de $X$-sch\'emas en groupes lisses :

\[\xymatrix{
       1 \ar[r] & \mu_n \ar[rr]  && SL_n \ar[rr]    &&PGL_n \ar[r]   & 1 \\
       1 \ar[r] &  \mu_n  \ar[u]^{ x \mapsto x} \ar[rr]^{x \mapsto
(x,x^{-1})}  \ar[d]_-{x \mapsto x^{-1}}  &&
       SL_n \times \G_m \ar[u] _{{\rm pr_1}} \ar[rr]^-{(u,v) \mapsto uv}
\ar[d]^{{\rm pr_2}} && GL_n  \ar[r]\ar[d]^{\rm det} \ar[u] & 1 \\
1\ar[r] & \mu_n \ar[rr]^{} &&\G_m \ar[rr]^-{x \mapsto x^n}  &&   \G_m \ar[r]  &
1.
\\
 }
\]

\subsection{Transformations \'el\'ementaires d'alg\`ebres d'Azumaya sur
une surface}

Soient $X$ un sch\'ema, $A \in \Az_n(X)$  et
$D \subset X$ un diviseur de Cartier effectif, $I_D \subset \O_X$
l'id\'eal inversible le d\'efinissant.
Supposons que la restriction
$A_D $ de $A$ \`a $D$ s'\'ecrive ${\End}(\V)$,
avec $\V$ un fibr\'e vectoriel sur $D$, et que l'on
dispose  d'un sous-fibr\'e vectoriel  $\F
\subset \V$ sur $D$ localement facteur direct. On construit
alors une autre alg\`ebre d'Azumaya $A' \in \Az_n(X)$, appel\'ee transform\'ee
\'el\'ementaire de
$A$ par rapport \`a $\F \subset \V$, de la fa\c con suivante.

 On consid\`ere  la sous-alg\`ebre $B $ des sections de $A$
qui pr\'eservent la filtration $\F \subset \V$.
Notons $i : D \hookrightarrow X$ l'immersion ferm\'ee naturelle. On a
la suite
exacte de $\O_X$-modules coh\'erents
$$ 0 \to B \to A \to i_*\Hom_D(\F,\Q) \to 0, \hskip2cm (3.2.1)$$
o\`u $ \Q = \V/\F$.
Notons que
$B$ contient $\O_X \subset A$.

Pour tout point $x$ de $X$, il existe un sch\'ema affine $U$
et un morphisme \'etale $U \to X$ d'image contenant $x$
tel qu'apr\`es restriction \`a $U$  on puisse \'ecrire
 $A=\End(V)$ avec $V_D \simeq \V$ et que de plus
il existe un scindage $V = F \oplus Q $ avec
$F_D=\F$.
Sur $U$, $A_U$ se lit
$$
\begin{pmatrix} \End(F) & \Hom(Q,F) \\ \Hom(F,Q)  & \End(Q) \end{pmatrix}
$$

On d\'efinit une alg\`ebre d'Azumaya $A'_U$ sur $U$ par
$$\begin{pmatrix} \End(F)  & I_D^{-1}\otimes \Hom(Q,F)  \\ I_D \otimes \Hom(F,Q)
  & \End(Q) \end{pmatrix}$$

On v\'erifie
que les alg\`ebres $A'_U$ pour divers $U$ se recollent
et d\'efinissent une alg\`ebre d'Azu\-maya $A' \in \Az_n(X)$, que l'on appelle
la transform\'ee \'el\'ementaire de $A$ le long de  $\F \subset \V$.
On v\'erifie sur la description locale ci-dessus
que l'on a la suite exacte de $\O_X$-modules coh\'erents
$$ 0 \to B \to A' \to i_*\Hom_D((I_D/I_D^2)\otimes \Q, \F) \to 0.\hskip2cm
(3.2.2)$$
Ici encore, l'inclusion $B \subset A'$ induit l'identit\'e de
$\O_X\subset B$  vers $\O_X \subset A'$.

(C'est une variante de constructions que l'on trouve dans d'autres
contextes :
transformations \'el\'ementaires entre ordres maximaux d'une alg\`ebre
simple centrale sur le corps des fractions d'un anneau de valuation
discr\`ete
\cite{[R]}, transformations \'el\'ementaires sur les fibr\'es vectoriels
(\cite{[HL]}).)

  \begin{lemm}[]\label{2.2.1} {\it 
Soient $X$ un sch\'ema  connexe, $D$ un
diviseur de Cartier effectif sur
$X$, et $i : D \subset X$ l'inclusion naturelle. Supposons donn\'ee une suite
exacte de $\O_X$-modules coh\'erents
$$ 0 \to V' \to V \to i_*Q \to 0$$
o\`u $V'$ et $V$ sont des fibr\'es vectoriels sur $X$ et $Q$ est un
fibr\'e vectoriel de rang $s$ sur $D$. L'application
d\'eterminant : $\det(V') \to \det(V)$ identifie le fibr\'e inversible
$\det(V')$ avec le sous-fibr\'e inversible $\det(V) \otimes \O_X(-sD)
\subset \det(V)$.}
\end{lemm}
 \noindent{\sc Preuve} --- On dispose de l'inclusion naturelle de fibr\'es
inversibles
$\det(V') \to \det(V)$. Pour \'etablir l'\'enonc\'e,
on peut supposer le sch\'ema $X$ local, soit $X=\Spec(R)$ et $D$
d\'efini par un \'el\'ement non inversible $\pi  \in R$.
On a une surjection de $R/\pi$-modules libres $V/\pi  \to Q$.
En scindant celle-ci on trouve une base ${\overline e}_1,\dots,{\overline e}_n$
du $R$-module libre $V/\pi$   telle que  ${\overline e}_{s+1},\dots,{\overline e}_n$
soit une base du noyau. Soient $(e_1, \dots, e_n)$ des relev\'es de ${\overline e}_1,\dots,{\overline
e}_n$ dans $V$. Alors $(e_1, \dots, e_n)$ est une base du $R$-module libre $V$ et
$(\pi e_1,\dots,\pi  e_s, e_{s+1},\dots, e_n)$ une base du $R$-module libre $V'$.
L'application $\Lambda^n V' \to \Lambda^nV$ envoie le  g\'en\'erateur
\'evident de $\Lambda^n V'$ sur $\pi^s.(e_1\wedge \dots \wedge e_n)$ qui est une base
de $\det(V) \otimes \O_X(-sD)$.
\cqfd

 \begin{prop}[]\label{2.2.2} Soient $X$ un sch\'ema  connexe, $A \in
\Az_n(X)$, avec $n$ inversible sur $X$, 
 et $D \subset X$ un diviseur de Cartier effectif. Supposons qu'il
existe un fibr\'e vectoriel $\V$ sur $D$ tel que
 $A_D = {\End}(\V)$, et supposons donn\'e 
 un sous-fibr\'e
vectoriel  $\F
\subset \V$ sur $D$ localement facteur direct, de rang constant $r$. Soit
$A'$ le transform\'e \'el\'ementaire de $A$ le long de
$\F
\subset
\V$.

 On a alors
$$ \cl(A') = \cl(A) -  r. [D] \in H^2(X,\mu_n),$$
o\`u $[D]$ d\'esigne l'image de  $D$
par l'application compos\'ee  $\Div(X) \to \Pic(X)/n \hookrightarrow
H^2(X,\mu_n)$.
\end{prop}

 \noindent{\sc Preuve} (P. Gille) ---    Soit $p: Y \to X$ le sch\'ema de
Severi-Brauer
associ\'e \`a $A$. Soit $i : D_Y=D\times_X Y \to Y$ l'immersion
naturelle.  Sur $Y$, il y a une suite exacte
$$ 0 \to V' \to V \to i_*\Q \to 0,$$
avec $V'$ et $V$ fibr\'es vectoriels sur $Y$ et $\Q$ le fibr\'e vectoriel
de rang $n-r$ sur $D_Y$ image r\'eciproque de $\V/\F$ par $p$.  On 
a $A_Y={\End}(V)$ et $A'_Y={\End}(V')$.
On a donc (Lemme \ref{2.1.1} (c)) $\cl(A_Y) = -\det(V) \in H^2(Y,\mu_n)$ et
$\cl(A'_Y)=-\det(V')
\in H^2(Y,\mu_n)$, o\`u l'on utilise tacitement l'inclusion $\Pic(Y)/n
\hookrightarrow
H^2(Y,\mu_n)$.
 D'apr\`es le lemme \ref{2.2.1}
 on a donc $\cl(A'_Y)=\cl(A_Y)  + (n-r). [D]_Y =\cl(A_Y)  -r [D]_Y\in
H^2(Y,\mu_n)$.
En analysant la suite spectrale de Leray pour le morphisme $p : Y \to X$
on montre que la fl\`eche naturelle $H^2(X,\mu_n) \to H^2(Y,\mu_n)$
est injective.  Ainsi $ \cl(A') = \cl(A) -  r. [D] \in H^2(X,\mu_n).$ \cqfd

On peut aussi \'etablir cette formule par un calcul local
via des cocycles de \v{C}ech (\cite{[dJ]}).

\begin{lemm}[]\label{2.2.4}  {\it
Soit $k$ un corps infini.  Soient $D/k$ une courbe
projective lisse g\'eom\'e\-tri\-quement connexe et $V$ un fibr\'e vectoriel de
rang $n$ sur $D$.  Supposons donn\'e  pour chaque point
$t$ dans un ensemble fini $T$ de points de $D(k)$ un
sous-espace vectoriel $F_{t }\subset V_{t}$ de dimension $1$.
Pour tout $m>0$ suffisamment grand il existe un
homomorphisme injectif
$ \O_D(-m) \to V$
dont l'image est localement facteur direct
dans $V$ et tel que l'image de $\O_D(-m)_t$ dans $V_t$
co\"{\i}ncide avec $F_t$. }
\end{lemm}
\noindent{\sc Preuve} --- Soit $W \subset V$ le sous-$\O_{D}$-module (localement libre)   de $V$
dont les sections apr\`es \'evaluation en $t$ donnent
des \'el\'ements de $F_{t}$. Pour $m \gg 0$, le fibr\'e vectoriel
$W(m)$ est engendr\'e par ses sections. Soit $E$ le $k$-vectoriel $H^0(D,W(m))$.
Notons encore $E$ l'espace affine d\'efini par $E$.
Soit $Z \subset D \times_{k}E$ le ferm\'e dont les points g\'eom\'etriques
sont les couples $(x,e)$ avec $e \in E$ tel que   $e_{x} = 0 \in W(m)_{x}$.
Comme $W$ est de rang au moins $2$,  
et que $W(m)$ est engendr\'e par ses sections globales,
pour $x$ point g\'eom\'etrique fix\'e,
l'ensemble des $e \in E$ satisfaisant $e_{x}=0$ est de codimension au
moins 2 dans $E$. La codimension de $Z$ dans $D \times_{k} E$ est
donc au moins 2, et le ferm\'e $Z_{1} \subset E$ qui est
l'adh\'erence de la projection de $Z$ dans $E$
est donc de codimension au moins 1. Choisissons un $k$-point de $E$ dans
le compl\'ementaire de $Z_{1}$. Ceci d\'efinit un homomorphisme 
$\O_{X} \to W(m)$ et donc un homomorphisme $\O_{X}(-m)  \to W \subset V$
satisfaisant les propri\'et\'es annonc\'ees. \cqfd

\subsection{ Rel\`evement des alg\`ebres}

 Soient $X$ un sch\'ema et $A \in \Az_n(X)$. La trace r\'eduite
$ \Tr  : A  \to \O_X$
est une application $\O_X$-lin\'eaire qui induit sur $\O_X \subset A$
la multiplication par $n$. On note $A^0 \subset A $ le noyau de la trace
r\'eduite.
Lorsque $n$ est inversible sur $X$, l'injection $\O_X \subset A$ 
est scind\'ee, et le quotient $A/\O_X$  est isomorphe \`a $A^0$.

  \begin{theo}[]\label{2.4.1}
Soit $X$ une surface connexe,  projective
et lisse
sur un corps  $k$
 alg\'ebri\-que\-ment clos. Soit $A \in \Az_n(X)$  avec
$n>1$ premier
\`a la caract\'eristique de $k$. Il existe une alg\`ebre d'Azumaya $A'
\in \Az_n(X)$, obtenue par transformation
\'el\'ementaire de $A$, telle que $\cl(A')=\cl(A) \in H^2(X,\mu_n)$
et que $H^2(X,A'/\O_X)=0$. 
\end{theo}

 \noindent{\sc Preuve} ---
Soit $\L$ un faisceau inversible sur $X$. L'inclusion $\O_X \subset A$
induit une inclusion $\L \subset A \otimes \L$ et donc une inclusion
$H^0(X,\L) \subset H^0(X,A \otimes \L)$.

\medskip
1)  Pour tout \'el\'ement $s$  de
$H^0(X,A\otimes \L)$ n'appartenant pas \`a $H^0(X,\L)$ il
existe une infinit\'e de $t \in X(k)$ tels que $s_t$ ne soit pas {\og scalaire \fg}, et
donc tel qu'apr\`es choix d'un isomorphisme $A_t \simeq End(V_t)$ 
il existe un sous-espace vectoriel $F_t \subset V_t$ de dimension 1
tel que $s(F_t)$ ne soit pas contenu dans $F_t \otimes \L_t$.
Comme $H^0(X,A\otimes \L)$ est de dimension finie, on en d\'eduit :

 Il existe un ensemble fini $T \subset X(k)$ et pour
chaque $t \in T$ un $k$-espace vectoriel $V_t$ de dimension $n$, un
sous-espace vectoriel $F_t \subset V_t$ de dimension $1$ et un
isomorphisme d'alg\`ebres $A_t \simeq End(V_t)$  tels que
l'inclusion $\L \subset A \otimes \L$ induise une \'egalit\'e
$$H^0(X,\L)=\{ s \in H^0(X,A \otimes \L), \hskip1mm s_t(F_t) \subset F_t \otimes\L_{t}
\hskip2mm {\rm pour \hskip1mm tout}
\hskip1mm t \} .$$

\medskip

2)  Soit $T \subset X(k)$ un ensemble fini de points.
Une variante du th\'eor\`eme de Bertini montre :
Il existe une courbe lisse connexe $D \subset X$  d'image
nulle dans $\Pic(X)/n\Pic(X)$ telle que $T \subset D$. 
(Il suffit de prendre une section convenable
 du faisceau $\O_X(nq)$
pour $q>0$ assez grand.)

\medskip

3) Choisissons $T$ comme en 1) et $D \subset X$ comme en 2).
D'apr\`es le th\'eor\`eme de Tsen, Ê{\it le groupe de Brauer 
du corps des fonctions d'une courbe d\'efinie sur un corps
alg\'ebriquement clos est nul. Le th\'eor\`eme \ref{injec} donne
alors $\br(D)=0$}\footnote{Dans la d\'emonstration du th\'eor\`eme de de Jong,
cet argument n'est pas utilis\'e, car on d\'eforme une alg\`ebre de classe triviale dans
le groupe de Brauer, voir le d\'ebut de la d\'emonstration de la proposition \ref{2.6.1}.}.
D'apr\`es le  lemme \ref{2.1.1} il existe donc
 un fibr\'e vectoriel
$\V$ de rang $n$ sur $D$  tel que $A_D \simeq \End_D(\V)$.
Pour $t \in T$, on a deux isomorphismes  $A_t \simeq End(V_t)$
et $A_{D,t} \simeq End_D(\V_t)$. Ces deux isomorphismes sont d\'eduits
l'un de l'autre via des isomorphismes $\eta_t : V_t \simeq \V_t$.
Pour tout $m>0$  assez grand, le lemme \ref{2.2.4}  
assure l'existence d'un   homomorphisme  $\varphi : \F \to \V$ ,
avec $\F=\O_D(-m)$,
 tel qu'en tout $t \in T$
l'image de $\varphi_t$ soit \'egale \`a $\eta_t(F_t)$. Soit $A' \in
\Az_n(X)$ le transform\'e \'el\'ementaire de $A$ le long de  $\varphi : \F
\to
\V$.

\medskip

4) Le choix de $D$, de classe nulle dans $\Pic(X)/n\Pic(X)$ et la
proposition \ref{2.2.2} impliquent
$$ \cl(A')=\cl(A) \in H^2(X,\mu_n).$$

\medskip

5)  La suite  exacte 
(3.2.1)  donne  apr\`es tensorisation avec $\L$
la suite exacte
$$ 0 \to B\otimes \L \to A\otimes \L  \to i_*\Hom_D(\F,\Q)\otimes \L_D \to 0 \hskip2cm
(3.4.1)$$ 
La fl\`eche compos\'ee
$$H^0(X,B\otimes \L) \to H^0(X,A\otimes \L)
\to H^0(D,i_*\Hom_D(\F,\Q) \otimes \L_D) \to \prod_{t \in T}
Hom(\F_t,\Q_t) \otimes \L_t$$
est nulle.
La premi\`ere fl\`eche est clairement injective. Le noyau de
la deuxi\`eme s'identifie par construction \`a $H^0(X,\L) \subset H^0(X,A\otimes \L)$.
On en conclut $H^0(X,\L)=H^0(X,B\otimes \L)$.
En tensorisant la suite exacte (3.2.2) par le fibr\'e inversible
$\L$ on obtient la suite exacte
$$ 0 \to B\otimes \L \to A'\otimes \L \to i_*\Hom_D((I_D/I_D^2)\otimes \Q, \F)\otimes
\L_D
\to 0.\hskip2cm (3.4.2)$$
On voit maintenant que le  noyau de
$$H^0(X, A'\otimes \L) \to
H^0(D, \Hom_D((I_D/I_D^2)\otimes \Q, \F) \otimes_D \L_D )$$
s'identifie \`a $H^0(X,\L) \subset H^0(X, A'\otimes \L)$.

Notons ${\mathcal N}_{X/D}$ le faisceau normal de $D$ dans $X$
et $E^*$ le  dual d'un $\O_D$-faisceau localement libre $E$.
Le $\O_D$-faisceau localement libre $\Hom_D((I_D/I_D^2)\otimes \Q, \F)
\otimes_D \L_D$ peut encore s'\'ecrire  ${\mathcal N}_{X/D} \otimes
\Q^*
\otimes
\F \otimes  \L_D$. Il s'injecte dans
le faisceau ${\mathcal N}_{X/D} \otimes \V^* \otimes
\F \otimes  \L_D$. On a donc aussi une injection
$$ H^0(D,\Hom_D((I_D/I_D^2)\otimes \Q, \F)
\otimes_D \L_D) \hookrightarrow H^0(D,{\mathcal N}_{X/D} \otimes \V^*
\otimes
  \L_D \otimes \F ).$$

 Observons alors qu'une fois fix\'es
$T, D, \V$ on peut choisir $m$  au point 3) aussi grand que l'on veut pour
d\'efinir $\varphi : \O_D(-m)=\F \subset \V$, et donc pour d\'efinir $A'$.
Si l'on choisit  $m$ tel que
 $$H^0(D, {\mathcal N}_{X/D} \otimes \V^* 
 \otimes  \L_D \otimes \O_D(-m))=0$$ et $\varphi$
et donc $A'$ avec un tel $m$, alors
$H^0(X,\L) = H^0(X, A' \otimes \L).$
Comme $n$ est inversible sur $X$, ceci implique
$H^0(X, A'^0\otimes \L)=0$.

\medskip

6)   Appliquons le r\'esultat obtenu
au faisceau inversible  canonique $\L=\omega_{X/k}=\Lambda^2(\Omega^1_X)$.
Cela donne une alg\`ebre d'Azumaya $A'\in \Az_n(X)$ qui est transform\'ee
\'el\'ementaire de $A$, qui satisfait $\cl(A)=\cl(A') \in H^2(X,\mu_n)$
et telle que $H^0(X,A'^0\otimes
\omega_{X/k})=0$.
L'application $A' \times A' \to \O_X$ d\'efinie par $(x,y) \mapsto
\Tr(xy)$ 
 induit  une dualit\'e parfaite entre les $\O_X$-modules $A'^0$
et $A'/\O_X$. 
Par dualit\'e de Serre sur la surface projective et lisse $X$,
le $k$-espace vectoriel $H^2(X,A'/\O_X)$ est dual de $H^0(X,A'^0\otimes
\omega_{X/k})$. Ainsi $H^2(X,A'/\O_X)=0$.
\cqfd

\bigskip

Le th\'eor\`eme \ref{2.4.1} va permettre d'utiliser les
th\'eor\`emes de rel\`evement suivants.

 \begin{prop}[ ]\label{2.4.2} Soit $X \to X'$ une
immersion ferm\'ee de sch\'emas d\'efinie par un id\'eal $I \subset \O_{X'}$ tel que $I^2=0$.
L'id\'eal $I$ est un  $\O_X$-module.
Soit $A$ une alg\`ebre d'Azumaya sur
$X$. Si l'on a $H^2(X, (A/\O_X) \otimes I)=0$, alors il existe une alg\`ebre
d'Azumaya $A'$ sur $X'$ induisant $A$ sur $X$. \cqfd
\end{prop}
ÊC'est un cas particulier d'un r\'esultat g\'en\'eral de Giraud 
(\cite{[Gi]},  Chap. VII, Th\'eor\`eme  1.3.1 et Remarque 1.3.1.2). \cqfd

\begin{theo}[]\label{2.4.3}
Soit $R$ un anneau local noeth\'erien complet
de corps r\'esiduel~$k$.
Soit $X/\Spec(R)$ un sch\'ema propre et plat.
 Soit $X_0/k$ la fibre sp\'eciale.
Si $A_0 \in \Az_n(X_0)$ satisfait  $H^2(X_0,A_0/\O_{X_0})=0$,
alors il existe une alg\`ebre d'Azumaya $A \in \Az_n(X)$ qui
induit $A_0$ sur $X_0$.
\end{theo}
Cela r\'esulte de la proposition \ref{2.4.2} et  du th\'eor\`eme
d'alg\'ebrisation des modules formels de Grothendieck (FGA, SGA 1, \cite{[EGA]} III.5
Th\'eor\`eme 5.1.4, \cite{[Gr]},  III, \S 3,  \cite{[I]}). \cqfd

 \begin{theo}[]\label{2.4.4}
Soit $R$ un anneau local hens\'elien
de corps r\'esiduel~$k$.
Soit $X/\Spec(R)$ un sch\'ema propre et plat.
 Soient  $X_0/k$ la fibre sp\'eciale et 
$A_0 \in \Az_n(X_0)$. Si l'on a $H^2(X_0,A_0/\O_{X_0})=0$, 
alors il existe une alg\`ebre d'Azumaya $A \in \Az_n(X)$ qui
induit $A_0$ sur $X_0$.
\end{theo}

C'est une cons\'equence formelle du th\'eor\`eme pr\'ec\'edent
et du th\'eor\`eme d'approxi\-ma\-tion
d'Artin \cite{[Ar1]}. Pour une situation similaire et quelques d\'etails
de plus, voir \cite{[CTOP]}, Theorem 1.8.
On commence par se r\'eduire au cas o\`u
$R$ est l'hens\'elis\'ee d'une $\Z$-alg\`ebre de type fini
en un id\'eal premier.
Le foncteur $S \to H^1(X \times_RS, PGL_n)$
de la cat\'egorie des $R$-alg\`ebres  commutatives vers la cat\'egorie
des ensembles est de pr\'esentation finie. Le th\'eor\`eme pr\'ec\'edent
donne une classe $\hat{\xi} \in H^1(X \times_R\hat{R}, PGL_n)$ d'image
$\xi_0 \in H^1(X_0,PGL_n)$. Le th\'eor\`eme d'Artin assure l'existence
de $\xi \in H^1(X,PGL_n)$ de m\^eme image que $\hat{\xi}$
dans $H^1(X_0,PGL_n)$, c'est-\`a-dire d'image $\xi_0$.
\cqfd

\subsection{Scindage des alg\`ebres d'Azumaya sur une surface et mise en
famille}

Soit $k$ un corps.
 Soit $X$ une $k$-surface projective, lisse, g\'eom\'e\-tri\-quement
int\`egre et $ A \in \Az_m(X)$ une alg\`ebre d'Azumaya sur $X$.
Soit $\L$ un faisceau inversible sur $X$. Notons  $L=\Spec(\Sym(\L^{-1}))  \to
X$ le fibr\'e en droites correspondant. \`A toute section $\sigma \in H^0(X,A
\otimes \L)$ on associe un id\'eal caract\'eristique $I_{\sigma}  \subset
\O_L$
 (\cite{[OP]}, Lemma 2.10).  Localement sur $X$ il est d\'efini de la fa\c con suivante.
 Soit $U=\Spec(R) $ un ouvert affine de $X$ sur lequel on dispose
 d'un isomorphisme $\O_U \simeq \L_U$. Soit $f \in \L(U)$
 l'image de $1$. On a alors   $L_U \simeq U[T]=\Spec(R[T])$ et l'id\'eal
$I_{\sigma}$
 est engendr\'e par le polyn\^ome caract\'eristique r\'eduit $P_{f,U}[T]$
 de $\sigma.f^{-1}$ qui est un \'el\'ement de
 l'alg\`ebre d'Azumaya  $A(U) \in \Az_m(U)$. Ce polyn\^ome est de degr\'e $n$ et
unitaire.
 Ainsi le ferm\'e $Y_{\sigma} \subset L$ d\'efini par $I_{\sigma}$
 est fini et plat sur $X$, de degr\'e $m$.

 \begin{theo}[M. Artin]\label{2.5.1} Soient $k$ un corps alg\'ebriquement clos, $n$ un entier
 premier \`a la carac\-t\'eristique de $k$, 
 $X$ une $k$-surface projective, lisse connexe,
$A \in \Az_m(X)$  telle que $A_{k(X)}$ est \`a division  et $\L$ un faisceau inversible.
Si  $\L$ est suffisamment ample et $\sigma$ est une section 
suffisamment g\'en\'erale du 
 fibr\'e vectoriel $A \otimes \L$,
alors
 la surface $Y_{\sigma}$ d\'efinie ci-dessus
est lisse et  connexe. Pour tout tel $\sigma$, l'image
r\'eciproque de
$A$ sur
$Y_{\sigma}$ a une classe triviale dans
$\br(Y_{\sigma})$. 
\end{theo}

Ce th\'eor\`eme est annonc\'e par de Jong (\cite{[dJ]},
\S 8). Une d\'emonstration en caract\'eristique nulle est donn\'ee par
Ojanguren et Parimala (\cite{[OP]}). Ceux-ci apportent de multiples
pr\'ecisions sur le type de ramification auquel on peut
de plus restreindre $Y_{\sigma} \to X$.
Je ne donnerai pas la d\'emonstration de ce th\'eor\`eme.
La premi\`ere assertion est un \'enonc\'e de type Bertini, qui utilise de
fa\c con essentielle le fait que la dimension de $X$ est $2$.
La seconde est une cons\'equence du th\'eor\`eme \ref{injec}.  \cqfd

\medskip

En utilisant ce th\'eor\`eme Ojanguren et Parimala montrent (\cite{[OP]})~:

 \begin{prop}[]\label{2.5.2} Soient $k$ un corps alg\'ebriquement clos de
caract\'eristique z\'ero, 
$X$ une $k$-surface projective, lisse, connexe,
$\alpha \in \Br(X)$.
Il existe une $k$-vari\'et\'e lisse connexe $W$  de dimension  $3$,
et des morphismes
\[\xymatrix{
W \ar[r]^{g} \ar[d]_f & X \\
\A^1&
}
\]
satisfaisant les propri\'et\'es suivantes :

{\rm (i)} Le morphisme $f$ est lisse \`a fibres connexes.

{\rm (ii)} Le morphisme $(g,f) : W \to X \times_k \A^1_k$ est quasi-fini et plat.

{\rm (iii)} Il existe un voisinage de $0 \in \A^1_k$ au-dessus duquel
$f$ est propre. En particulier la surface $Y=f^{-1}(0)$ est projective, lisse, connexe.

{\rm (iv)} On a
 $g^*(\alpha)= 0 \in \Br(Y)$.

{\rm (v)} La fibre $W_1=f^{-1}(1)$ est non vide, et la restriction
de $g : W \to X$ \`a $W_1$ est une immersion ouverte.
\end{prop}

  \noindent{\sc Preuve (esquisse)} ---  D'apr\`es le th\'eor\`eme \ref{purdim2}, il existe une alg\`ebre
  d'Azumaya $A \in \Az_{m}(X)$, avec $m$ convenable, de classe $\alpha \in \Br(X)$ et telle
  que $A_{k(X)}$ est \`a division.
  Soient $\L, L$ et $\sigma \in H^0(X,\L)$ comme
au th\'eor\`eme \ref{2.5.1}.  On choisit des sections globales distinctes $w_1,\dots,w_m$ de $L$.
On consid\`ere le fibr\'e en
droites
$L\times_k\A^1_k \to X \times_k\A^1_k$, o\`u $\A^1_k=\Spec(k[t])$.
Pour $U=\Spec(R) \subset X$ et $f \in \L(U)$ comme ci-dessus, on 
consid\`ere l'id\'eal de de $R[T,t]$ engendr\'e par
$$Q_{f,U}(t,T)=(1-t)P_{f,U}(T)+t(T-w_1/f)\dots(T-w_m/f).$$
On v\'erifie que ces diff\'erents id\'eaux se recollent en un id\'eal
sur le sch\'ema $L\times_k\A^1_k$. Soit $Z \subset L\times_k\A^1_k$
le ferm\'e d\'efini par cet id\'eal.
La projection $Z \to X\times_k\A^1_k$ est finie et plate de degr\'e $m$.
L'application compos\'ee $f : Z \to X\times_k\A^1_k \to \A^1_k$ est donc propre
et plate.
Sa fibre en $t=0$ est la surface projective, lisse, connexe $Y_{\sigma}$.
La fibre g\'en\'erique de $f$ est donc
 lisse
et g\'eom\'etriquement int\`egre.
Comme $Z \to X$ est plat,
ceci implique que $Z$ est int\`egre.
La fibre au-dessus de $t=1$ contient $n$ composantes ferm\'ees distinctes
chacune birationnellement isomorphe \`a $X$ via la projection sur $X$.
En enlevant un ferm\'e convenable dans $Z$, on obtient $W$
comme annonc\'e dans la proposition.
\cqfd
\begin{rema}\label{2.5.3}
 Dans \cite{[dJ]}, le th\'eor\`eme \ref{2.5.1} n'est pas utilis\'e.
L'id\'ee de base, qui est de d\'eformer un polyn\^ome en une variable, unitaire,
 de degr\'e $m$,
qui d\'efinit l'extension $k(Y)/k(X)$, en un polyn\^ome unitaire s\'eparable de m\^eme
degr\'e avec toutes ses racines dans $k(X)$, est la m\^eme,
mais la construction d'une bonne famille satisfaisant des propri\'et\'es
plus faibles que celles de la proposition \ref{2.5.2}, est alors
plus d\'elicate (\cite{[dJ]}, \S 4). De Jong y a en particulier recours
au {\og th\'eor\`eme de la fibre r\'eduite \fg}  (\cite{[BLR]}, \cite{[dJS2]}),
qui avait d\'ej\`a \'et\'e utilis\'e dans \cite{[dJS1]}.
\end{rema}

\subsection{ Le th\'eor\`eme de de Jong dans le cas non ramifi\'e}

 \begin{prop}[]\label{2.6.1}
Soit $X$ une surface projective, lisse, connexe sur un corps $k$
 alg\'e\-bri\-quement clos, $n>0$ entier inversible dans $k$.
Soit $\alpha \in \Br(X) \subset \Br(k(X))$, d'exposant $n$.

Supposons qu'il existe une vari\'et\'e int\`egre $W$ de dimension 3
et des morphismes
\[\xymatrix{
W \ar[r]^{g} \ar[d]_f & X \\
\A^1&
}
\]
satisfaisant les  propri\'et\'es suivantes :

{\rm (i)} Le morphisme $f$ est lisse \`a fibres 
connexes.

{\rm (ii)} Le morphisme $f$ est propre au voisinage du point $0 \in \A^1$,
de fibre $W_0=f^{-1}(0)$, et $g$ induit un morphisme fini et plat  $g_{W_0} : W_0
\to X$
de surfaces lisses.

{\rm (iii)} La restriction de $g$ \`a la fibre $W_1=f^{-1}(1)$ est une
immersion ouverte $W_1 \hookrightarrow X$.

{\rm (iv)} $g^*(\alpha)_{W_0} = 0 \in \Br(W_0)$.

Alors l'indice de $\alpha_{k(X)} \in \Br(k(X))$ est \'egal \`a $n$.
\end{prop}

 \noindent{\sc Preuve} --- Soit $\eta \in H^2(X,\mu_n)$ d'image
$\alpha \in {}_{n}\Br(X)$.
Notons $Y=W_0$.
Soit $\L_Y$ un faisceau inversible sur $Y$ dont la classe
par l'application compos\'ee $\delta : \Pic(Y)/n \to H^2(Y,\mu_n)$
est l'oppos\'e de $g^*(\eta)_Y$. Qu'un tel faisceau inversible
existe r\'esulte de la suite exacte de Kummer et de l'hypoth\`ese
(iv). Consid\'erons l'alg\`ebre d'Azumaya   $$B_0={ \End}(\L_Y
\oplus (\O_Y)^{n-1}).$$ {\it Elle est de degr\'e $n$.} 
 D'apr\`es le lemme \ref{2.1.1} (c) on a
$\cl(B_0)=-\delta(\L_Y)=g^*(\eta)_Y$. D'apr\`es le th\'eor\`eme \ref{2.4.1}, 
une transformation \'el\'ementaire convenable de $B_{0}$ produit
 une alg\`ebre $A_{0} \in Az_{n}(Y)$ satisfaisant  $\cl(A_0)=\cl(B_0)=g^*(\eta)_Y \in 
H^2(Y,\mu_n)$ et $H^2(Y,A_0/\O_Y)=0$.
{\it C'est un point-cl\'e de la d\'emonstration}.
Cette nullit\'e permet de d\'eformer l'alg\`ebre $A_0 \in \Az_n(Y)$. En utilisant
le  {\it th\'eor\`eme  d'alg\'ebrisation}  \ref{2.4.3} puis le
{\it th\'eor\`eme d'approximation}  \ref{2.4.4},
 on voit qu'il existe un voisinage \'etale
connexe
$h : (C',0) \to (C,0)$ du point $0 \in C=\A^1_k$, tel que
$W'=W\times_CC' \to C'$ soit propre et lisse, et une alg\`ebre d'Azumaya
$A \in \Az_n(W')$ telle que la fibre de $A$ au-dessus de  $0 \in
C'$ soit isomorphe \`a l'alg\`ebre $A_0$ sur $Y$.
(Notons que l'alg\`ebre $A$ n'a pas de raison d'avoir une classe
triviale dans le groupe de Brauer de $W'$, m\^eme
apr\`es compl\'etion en $0 \in C'$.)

La r\'eduction de $\cl(A) \in H^2(W',\mu_n)$
au-dessus de $0 \in C'$
s'identifie \`a $\cl(A_0)=g^*(\eta)_Y$, c'est-\`a-dire
\`a la r\'eduction de l'image au-dessus de $0 \in C'$
de la classe globale obtenue \`a partir de $\eta  \in H^2(X,\mu_n)$
en prenant l'image r\'eciproque par le morphisme compos\'e $g' : W'\to W \to X$.

En appliquant  le {\it th\'eor\`eme de changement de base propre en cohomologie \'etale}
(\cite{[Mi]}, VI.2.7) au morphisme
$W' \to C'$, on voit qu'en remplac\c ant $(C',0)$ par un autre
voisinage \'etale, encore not\'e $(C',0)$, on peut
assurer  $\cl(A)=g'^*(\eta) \in H^2(W',\mu_n)$.

On peut \'etendre le morphisme $C' \to C$ en un morphisme
{\it fini} de courbes lisses $D \to C=\A^1_k$. Le morphisme $W\times_CD \to D$
est lisse, et {\it comme $k$ est alg\'ebriquement clos} il existe un point
$1 \in D$ au-dessus de $1$, la fibre $W'_1$ en ce point est lisse, connexe,
et la projection compos\'ee $W\times_CD \to W \to X$ induit sur
 $W'_1$ une immersion ouverte $W'_1 \subset X$.

Soit $R$ l'anneau local de  $W\times_CD$ au point g\'en\'erique de $W'_1$.
C'est un anneau de valuation discr\`ete, de corps des fractions $K$
le corps des fonctions de $W'$, de corps r\'esiduel isomorphe
\`a $k(X)$. De plus, la projection $W\times_CD \to W \to X$
induit une inclusion $k(X) \hookrightarrow R$ qui compos\'ee avec l'application
de r\'eduction modulo l'id\'eal maximal de $R$, soit $R \to k(X)$,  est
l'identit\'e de $k(X)$. Soit $A_K \in \Az_n(K)$  l'image de
de   $A \in \Az_n(W')$. L'\'egalit\'e  $\cl(A)=g'^*(\eta) \in H^2(W',\mu_n)$
implique que la classe $[A_{K}] \in \Br(K)$ 
co\"{\i}ncide avec l'image de $\alpha_{k(X)} \in \Br(k(X))$
via la fl\`eche compos\'ee $k(X) \to R \to K$. En particulier,
$[A_K]$ appartient \`a $\Br(R) \subset \Br(K)$.
L'alg\`ebre $A_K$ est de degr\'e $n$.  La combinaison de ces deux derniers
faits  et du th\'eor\`eme \ref{purdim2} (dans le cas particulier d'un anneau de valuation discr\`ete)
 assure que la sp\'ecialisation $\Br(R) \to \Br(k(X))$ envoie
$A_K$ sur la classe d'une alg\`ebre d'indice divisant $n$. Comme cette sp\'ecialisation
est  \'egale \`a $\alpha_{k(X)} $, et que cette derni\`ere est d'exposant $n$,
on voit que l'indice de  $\alpha_{k(X)} $ est $n$. \cqfd

Combinant les propositions \ref{2.5.2} et \ref{2.6.1},
on obtient une d\'emonstration en caract\'eristique nulle du th\'eor\`eme~:

 \begin{theo}[de Jong]\label{2.6.3}
Soit $k$ un corps alg\'ebriquement clos. Soit
$X/k$ une surface connexe, projective et lisse sur $k$. Soit $\alpha \in
\Br(X)$ d'exposant $n>0$ premier \`a la caract\'eristique de $k$. Alors
l'indice de $\alpha_{k(X)}$ est \'egal \`a $n$, et $\alpha$ est
repr\'esent\'e par une alg\`ebre d'Azumaya $A$ sur $X$ de degr\'e $n$. 
\end{theo}

La derni\`ere assertion est une cons\'equence du reste du th\'eor\`eme
et du th\'eor\`eme \ref{purdim2}.
 Une d\'emonstration de la proposition
 \ref{2.5.2} en caract\'eristique positive premi\`ere \`a l'exposant
 de l'alg\`ebre $A$ permettrait d'\'etendre la d\'emonstration
 ici d\'ecrite. De Jong quant \`a lui \'etablit le r\'esultat dans cette
g\'en\'eralit\'e en utilisant la mise en famille \'evoqu\' ee \`a la
remarque \ref{2.5.3} ci-dessus.

\subsection{Le cas ramifi\'e r\'esulte du cas non ramifi\'e}

 \begin{theo}[de Jong]\label{2.7.1}
Soit $k$ un corps alg\'ebriquement clos. Soit
$X/k$ une surface connexe, projective et lisse sur $k$. Soit $\alpha \in
\Br(k(X))$ d'exposant $n>0$ premier \`a la caract\'eristique de $k$.
Alors $\ind_{k(X)}(\alpha)=n.$
\end{theo}

  \noindent{\sc Preuve} --- La d\'ecomposition des corps gauches en produits
tensoriels de corps gauches d'indices premiers entre eux permet de se ramener
au cas o\`u $n=l^r$ avec $l$ premier.  
Montrons  le r\'esultat par r\'ecurrence sur $r$, en supposant le r\'esultat
connu pour $r=1$.  L'exposant de $\beta=l.\alpha \in \Br(k(X))$ 
est  $l^{r-1}$. Il existe donc une extension de corps $E/k(X)$
de degr\'e $l^{r-1}$ telle que $l.\alpha_E=\beta_E=0 \in \Br(E)$. La r\'esolution
des singularit\'es des surfaces alg\'ebriques donne l'existence d'une
surface projective lisse connexe $Y$ sur $k$ telle que $E=k(Y)$.
Le r\'esultat pour $n=l$  premier appliqu\'e \`a la surface $Y$
assure l'existence d'une extension de corps $F/E$ de degr\'e $l$
telle que $(\alpha_E)_F=0 \in \Br(F)$. Ainsi $\alpha_F=0$, et l'extension
de corps $F/k(X)$ est de degr\'e $l^r$. Cette r\'eduction  
est classique  (\cite{[A]}, p. 175; \cite{[CTG]}, p.
132).

On suppose  d\'esormais  $n=l$
premier, distinct de la caract\'eristique de $k$.

On va construire des morphismes
\[\xymatrix{
W \ar[r]^{g} \ar[d]_f & X \\
\A^1&
}
\]
de vari\'et\'es 
lisses connexes
satisfaisant les propri\'et\'es suivantes :

(i) Le morphisme $f : W \to \A^1_k$ est lisse.

(ii) La fibre g\'en\'erique g\'eom\'etrique $W_{\overline{\eta}}$
de $f : W \to \A^1_k$  est projective et connexe.

(iii) La fibre  $W_0 \subset W$
de $f$ en $0 \in \A^1_k$ est non vide et la restriction de $g : W \to X$
\`a   $W_0 \subset W$   est un  morphisme birationnel $W_0 \to X$.

(iv) L'image r\'eciproque $r^*(\alpha)$  de $\alpha$ par l'application
compos\'ee
$r : W_{\overline{\eta}} \to W_{\eta} \to W \to X$, o\`u la derni\`ere
fl\`eche est $g$, est non ramifi\'ee sur la surface
$W_{\overline{\eta}}$  (on note  $\eta$ le point g\'en\'erique de $\A^1_k$
et ${\overline{\eta}}$ un point g\'en\'erique g\'eom\'etrique).

\medskip

Supposons ces morphismes construits. De la propri\'et\'e (iv) et du th\'eor\`eme
\ref{2.6.3}  on d\'eduit qu'il existe  une extension finie de corps 
$K=k(\A^1)  \hookrightarrow
L $ telle que la restriction de $\alpha$ \`a $W_{\eta}\times_KL$
soit   de degr\'e $l$. Le corps $L$
est le corps des fonctions d'une courbe affine lisse $C$ munie
d'un $k$-morphisme fini $C \to \A^1_k$. {\it Comme $k$
est alg\'ebriquement clos}, il existe un $k$-point $M \in C$
au-dessus de $0 \in \A^1_k$.
La vari\'et\'e $Z=W\times_{\A^1}C$ est lisse, la fibre g\'en\'erique
de $Z \to C$ est g\'eom\'etriquement int\`egre, la vari\'et\'e
$Z$ est donc connexe.
La fibre de $Z_M$ de $Z \to C$ au-dessus du point $M$
est lisse, connexe,   la restriction de la projection
$Z \to X$ \`a la fibre $Z_M$ est un morphisme birationnel $Z_M \to X$.

Soit $R$ l'anneau local de $Z$ au point g\'en\'erique $\xi$ de
$Z_M$. Le corps des fractions de $R$ est le corps des fonctions
de $Z$, le corps r\'esiduel $\kappa$ est isomorphe \`a $k(X)$, plus
pr\'ecis\'ement le compos\'e $\Spec (\kappa) \to \Spec (R)
\to W \to X$ s'identifie \`a l'inclusion du point g\'en\'erique
de $X$ dans $X$. L'image r\'eciproque sur le corps des fonctions
de $Z$ de la classe $\alpha \in \Br(k(X))$, via la projection
$Z \to X$, est une classe dans $\Br ( k(Z))$ non ramifi\'ee au point
g\'en\'erique $\xi$  de $Z_M$, dont la r\'eduction sur $\kappa$
s'identifie
\`a $\alpha$.
Pour \'etablir le th\'eor\`eme \ref{2.7.1}, il suffit alors d'appliquer le th\'eor\`eme
\ref{purdim2} dans le cas particulier  du spectre d'un anneau 
de valuation discr\`ete.

\medskip

Il reste donc \`a construire un morphisme  $h=(g,f) : W \to X \times \A^1_k$
satisfaisant  les propri\'et\'es (i) \`a (iv).
Quitte \`a remplacer la surface $X$ par un \'eclat\'e, on peut
supposer (\cite{[Li]})
que le lieu de ramification de la classe $\alpha \in \Br(k(X))$
est un diviseur $D$ strictement \`a croisements normaux
(composantes lisses, par tout point il passe au plus deux composantes,
et dans ce cas elles se coupent transversalement).

Soit $\A^1=\Spec (k[t])$.
On cherche un rev\^etement g\'en\'eriquement cyclique
 $h=(g,f) : W \to X \times \A^1_k$
de degr\'e $l$, donn\'e par la racine $l$-i\`eme
d'une fonction rationnelle $f_t$ sur $X \times  \A^1$,
de telle sorte qu'au-dessus
du point g\'en\'erique g\'eom\'etrique de $\A^1$,
on ait d\'etruit la ramification de $\alpha$, et de
telle sorte qu'au-dessus du point $0 \in \A^1_k$
la fibre  de $f$ contienne une composante de multiplicit\'e 1
birationnelle \`a $X$. En bref, on veut que la sp\'ecialisation $f_0$
de $f$ soit une puissance  $l$-i\`eme, et on veut que
le diviseur de $f$ contienne le diviseur $D$ \`a l'ordre 1 (dans le cas mod\'er\'e, {\og la
ramification avale la ramification \fg} -- lemme d'Abhyankar).

\begin{lemm}[]\label{2.7.3}
{\it Etant donn\'e un diviseur  $D \subset X$ strictement \`a
croisements
normaux, il existe des diviseurs $E$ et $E'$ effectifs lisses tels  que
$D+E$ soit lin\'eairement \'equivalent \`a $l(D+E')$,
que deux de $D$,$E$ et $E'$ n'aient pas de composante commune,
et que $D+E+E'$ soit \`a croisements normaux.}
\end{lemm}

 \noindent{\sc Preuve} ---  Soit ${ M}$ un fibr\'e inversible ample sur $X$.
Pour $r$ entier suffisamment grand, le faisceau
$\O_X((l-1)D) \otimes { M}^{\otimes rl}$ est tr\`es ample,
il existe une section de ce faisceau qui d\'efinit sur $X$ un diviseur
lisse $E$ qui coupe $D$  transversalement (noter que les supports
de $D$ et de $E$ se rencontrent n\'ecessairement, le diviseur $D+E$
n'est pas lisse).
En outre il existe une section de $M^{\otimes r}$ dont le lieu des z\'eros
est un diviseur $E'$ lisse transverse \`a $D+E$. \cqfd

\medskip

 Soit ${ L}=\O_X(-D-E')$.
Soit $s_0$ une section de ${ L}^{\otimes {-l}}$ de lieu des z\'eros
$D+E$.
Soit $s_1$ une section de ${ L}^{-1}$ de lieu des z\'eros $D+E'$.
On dispose alors de la section $s_1^l = s_1^{\otimes l}$
de $L^{\otimes {-l}}$.

 Soit $\tilde{ L}$ l'image inverse de ${ L}$
sur $Y=X \times \A^1$. Soit $s_t$ la section de
 $\tilde{ L}^{\otimes {-l}}$ d\'efinie par
$$s_t = t^ls_1^l + (1-t)^ls_0,$$
o\`u, par abus de notation, $s_0$ et $s_1$ sont les images r\'eciproques
sur $X \times \A_1$ de $s_0$ et $s_1$.

 On d\'efinit alors, de fa\c con classique, un rev\^etement
g\'en\'eriquement cyclique $Y$ de $ X \times \A^1$ de la fa\c con
suivante. On d\'efinit une structure de $\O_{X \times \A^1}$-alg\`ebre sur
le fibr\'e vectoriel
$$\tilde{N}= \oplus_{i=0}^{i=l-1} \tilde{ L}^{\otimes i}$$
en utilisant la fl\`eche $\tilde{ L}^{\otimes l} \to \O_{X \times
\A^1}$ donn\'ee par  la section $s_t$ de $\tilde{ L}^{\otimes {-l}}$.

 Soit $U=\Spec (A)$ un ouvert affine de $X$ sur lequel ${ L}$ admet une
trivialisation ${ L}\simeq \O_U$. Au-dessus de $U$, le
rev\^etement d\'ecrit ci-dessus est isomorphe \`a
$$\Spec (A[t][x]/(x^l-(t^lf_1^l+ (1-t)^lf_0))),$$
pour $f_0$ et $f_1$ dans $A$
convenables.

On v\'erifie alors facilement :

(a) La vari\'et\'e  $Y$ est int\`egre.

(b) La fibre g\'en\'erique $Y_{\eta}$ du morphisme compos\'e $Y \to X
\times \A^1 \to
\A^1$ est g\'eom\'etriquement int\`egre.

(c) La fibre $Y_1$ du morphisme $Y \to X \times \A^1 \to \A^1$
au-dessus du point
$t=1$ a $l$ composantes, chacune
de multiplicit\'e 1, chacune isomorphe \`a $X$, l'isomorphisme
\'etant induit par le morphisme compos\'e $Y_1 \subset Y \to X$.

La $K$-vari\'et\'e $Y_{\eta}$ ne saurait \^etre lisse,
car le lieu des z\'eros de la section $s_t$ n'est pas lisse :
il est somme de $D\times_kK$ et d'un diviseur effectif non nul $E_t$.
Le lieu des z\'eros  de $s_0$ sur $Y_0=X$ est $D+E$, qui est
strictement \`a croisements normaux. Pour $t$ g\'en\'eral, ceci implique
(Bertini) la m\^eme propri\'et\'e pour la d\'ecomposition $D+E_t$
du lieu des z\'eros de $s_t$.
Ceci assure que  pour tout tel $t$, la vari\'et\'e $Y_t$
est normale et n'a que des singularit\'es de type $A_{l-1}$ : une
\'equation locale en un point singulier est analytiquement du type
$A[x]/x^l-uv$, o\`u $u,v$ sont des \'el\'ements d'un syst\`eme r\'egulier
de param\`etres.
On consid\`ere alors la r\'esolution minimale des singularit\'es
de la fibre g\'en\'erique $Y_{\eta}/K$, soit $Z_{\eta}/K$.
Cette r\'esolution est lisse.
La r\'esolution $Z_{\eta} \to Y_{\eta}$ est donn\'ee par l'\'eclatement
d'un certain id\'eal coh\'erent $I_{\eta}$ sur $Y_{\eta}$ de support les
points singuliers. On peut trouver un id\'eal coh\'erent $I$ sur $Y$
dont le lieu des z\'eros est de dimension au plus 1 et dont la restriction
\`a la fibre g\'en\'erique est $I_{\eta}$. Soit $Z \to Y$ l'\'eclat\'e de
$Y$ au moyen de l'id\'eal $I$.
Il existe alors un ouvert $W$ de $Z$ qui satisfait
les propri\'et\'es (i), (ii) et (iii) de l'\'enonc\'e
(pour la propri\'et\'e (iii), il suffit d'observer
qu'au-dessus des points g\'en\'eriques des composantes
des fibres de $Y \to \A^1$, l'\'eclatement via $I$ ne modifie rien;
on utilise alors le r\'esultat pour $Y \to \A^1$).
Il reste \`a \'etablir le point (iv).

La technique ici est standard 
(\cite{[Ar3]}, \cite{[FS]}, \cite{[CTOP]}).
Revenons \`a la situation envisag\'ee dans la d\'emonstration
principale.
Soit $L$ une cl\^oture alg\'ebrique de $K=k(t)$.
Le morphisme  de $L$-surfaces projectives, lisses, connexes
  $Z_L \to X_L$ est donn\'e au point g\'en\'erique de $X_L$ 
par l'extraction de la racine $l$-i\`eme de la fonction
rationnelle $f_t=s_t/s_1^l$, dont le diviseur est
$D+E_t-l(D+E')$.
Le  lieu de ramification de $\alpha_K \in \Br(X_K)$
est contenu dans $D_K$. Le diviseur $D+E_t$ est r\'eduit, strictement \`a
croisements normaux.

Soit $M$ un point de codimension $1$ de $Z_L$. Si son image
par $Z_L \to X_L$ n'est pas contenue dans $D_L$, alors
$\alpha_{L(Z)}$ est non ramifi\'e en $M$.
Si l'image de $M$ est un point g\'en\'erique $N$ d'une composante
de $D_K$, alors $\alpha_{L(Z)}$ est non ramifi\'e en $M$.
On a en effet le diagramme commutatif
\[\xymatrix{
\Br(L(Z)) \ar[r]^{}  & H^1(\kappa_M,{\bf Q}/\Z)   \\
\Br(L(X)) \ar[r]^{} \ar[u]_{}& H^1(\kappa_N,{\bf Q}/\Z) \ar[u]^{}  
 }
\]
o\`u la fl\`eche $H^1(\kappa_N,{\bf Q}/\Z) \to H^1(\kappa_M,{\bf Q}/\Z)$
est la multiplication par l'indice de ramification
en $N$, qui est $l$.

Supposons maintenant que l'image de $M$ soit un point ferm\'e $N$ de
$X_L$.
Si le point  $N$ n'appartient pas \`a $D $, alors
$\alpha_{L(X)}$ s'\'etend en une alg\`ebre d'Azumaya au voisinage
de $N$
donc  $\alpha_{L(Z)}$ est non ramifi\'e en $M$.
Supposons que le point $N$ appartienne \`a $D$ mais soit sur  une
unique composante de $D$. Alors (proposition \ref{1.3.2}) il existe un voisinage
affine $\Spec (R)$ de
$N$ dans $X_L$,
$u \in R^*$ et $s \in R$ tels que $s=0$ d\'efinisse exactement $D $ sur
$\Spec (R)$ et   que $\alpha_{L(X)} - (u,s)$
appartienne \`a $\Br(R)$ (on omet
ici 
 l'indice $\zeta_l$ dans la notation
d'une alg\`ebre cyclique).
On a donc $\partial_M(\alpha_{L(Z)})=\partial_M((u,s))={\overline
u}^{v_M(s)} \in \kappa_M/\kappa_M^{*l}$. La classe de $u$ dans
$\kappa_M/\kappa_M^{*l}$ est l'image de la classe de $u$ dans
$\O_{X_L,N}^*/\O_{X_L,N}^{*l}$, et l'application naturelle
  $\O_{X_L,N}^*/\O_{X_L,N}^{*l}\to \kappa_M/\kappa_M^{*l}$
se factorise par $\kappa_N^*/\kappa_N^{*l}=1$ (le corps
r\'esiduel $\kappa_N$ est s\'eparablement clos de caract\'eristique
diff\'erente de $l$.)
Supposons maintenant que le point $N$ soit l'intersection de deux
composantes de $D$. Il n'est donc pas sur $E_t$.
Il existe un voisinage affine $\Spec (R)$ de $N$ dans $X_L$,
$s,t \in R$ engendrant l'id\'eal maximal de $N$,
tels que $st=0$ d\'efinisse exactement $D$ sur $\Spec (R)$
et que $E_t$ ne rencontre pas $\Spec (R)$. Sur $\Spec (R)$
on a donc  $f_t=csth^l \in L(X)^*$ avec $c \in R^*$ et $h \in L(X)^*$.

D'apr\`es la proposition  \ref{1.3.2},
quitte  \`a restreindre
$\Spec (R)$ il existe  des unit\'es $u,v \in R^*$ et $r\in \Z$
tels  que
$\alpha+(s,u)+(t,v)+r(s,t)$ appartienne \`a $\Br(R)$.
On a
$$(s,t)= (s,f_t c^{-1}s^{-1}h^{-l})= (s,f_t)
+(s,-c^{-1}) \in \Br(L(X)),$$
o\`u l'on a utilis\'e la formule $(s,-s)=0$.
On peut donc \'ecrire $\alpha \in \Br(L(X))$
comme la somme d'\'el\'ements du type :
non ramifi\'e, $(s,f_t)$, $(s,u)$ et $(t,v)$ avec
$u$ et $v$ unit\'es. On a $L(Z)=L(X)(f_t^{1/l})$,
donc $(s,f_t)_{L(Z)}=0$. Le m\^eme argument que ci-dessus
montre que $(s,u)_{L(Z)}$ et $(t,v)_{L(Z)}$ ont des
r\'esidus triviaux en $M$. Ainsi $\alpha_{L(Z)}$ est
non ramifi\'e en $M$.

Ceci \'etablit le point (iv) et ach\`eve la d\'emonstration 
du th\'eor\`eme de de Jong. \cqfd

\bigskip

\begin{rema}\label{1-connexe}
 Le th\'eor\`eme de de Jong peut se reformuler ainsi :
 
 {\it Soit $K$ un corps de fonctions de deux variables sur un corps
 alg\'ebriquement clos. Soit  $A$ une alg\`ebre simple centrale sur $K$
 de degr\'e $nm$ et d'exposant $n$ premier \`a la caract\'eristique de $K$. Alors la $K$-vari\'et\'e de
 Severi-Brauer g\'en\'eralis\'ee $SB(A,n)$ poss\`ede un $K$-point.}
 
 On peut se demander ce qui fait que la $K$-vari\'et\'e $SB(A,n)$
 a automatiquement un  point rationnel. Dans un travail en pr\'eparation,
 de Jong et Starr \'etudient la notion de $1$-connexit\'e rationnelle et
 donnent des conditions suffisantes pour qu'une $K$-vari\'et\'e
 projective et lisse poss\`ede automatiquement un point rationnel
 sur $K$, corps de fonctions de deux variables sur un corps alg\'ebriquement clos.
 Leur r\'esultat est un analogue du r\'esultat de Graber, Harris, Starr, de Jong sur
 l'existence de points rationnels sur les vari\'et\'es (s\'eparablement) rationnellement
 connexes d\'efinies sur un corps de fonctions d'une variable. Les hypoth\`eses
 mises par de Jong et Starr dans le cas d'un  corps de fonctions de deux  variables
 sont assez contraignantes, mais elles s'appliquent aux vari\'et\'es $SB(A,n)$.
\end{rema}

\section{Cons\'equences pour les groupes alg\'ebriques lin\'eaires}

Soit $K$ un corps de fonctions de deux variables sur un corps
alg\'ebriquement clos de caract\'eristique z\'ero.
Pour un tel corps, les deux propri\'et\'es g\'en\'erales
suivantes sont donc satisfaites~:

\medskip

(1) {\it C'est un corps de dimension cohomologique $\cd(K) \leq 2$
(\cite{[S2]}).}

(2) {\it Sur tout corps extension finie de $K$, indice et exposant
des alg\`ebres simples centrales co\"{\i}ncident}.

\medskip

Ces   deux propri\'et\'es sont satisfaites par les corps mentionn\'es
au paragraphe 2.4, du moins par ceux qui ne sont pas formellement r\'eels :

\begin{itemize} 

\item[$\bullet$]  les corps $p$-adiques

\item[$\bullet$]  les corps de nombres totalement imaginaires

\item[$\bullet$]  les corps de la forme $F((t))$, o\`u    $F$ est un corps
de caract\'eristique  z\'ero  de dimension cohomologique $1$

\item[$\bullet$]  le corps des fractions d'un anneau
local int\`egre, excellent, hens\'elien, de dimension $2$, \`a corps r\'esiduel 
alg\'ebriquement clos de caract\'eristique z\'ero (\cite{[Ar3]},
\cite{[FS]}, \cite{[CTOP]}).
\end{itemize}

 \medskip

Un corps   de fonctions de deux variables sur un corps
alg\'ebriquement clos de carac\-t\'eristique z\'ero,
est un
corps $C'_2$ (th\'eor\`eme de Tsen-Lang), ce qui, d'apr\`es
Merkur'ev et Suslin
est une propri\'et\'e plus forte
que la propri\'et\'e (1). Comme mentionn\'e au \S 2, une autre application
du th\'eor\`eme de Merkur'ev et Suslin montre que la propri\'et\'e
$C'_2$ implique 
qu'exposant et indice
co\"{\i}ncident pour les alg\`ebres 2-primaires et 3-primaires.
Introduisons la condition suivante, plus faible que (2).
\medskip

(2') {\it Sur tout corps $L$ extension finie de $K$, indice et exposant
des $L$-alg\`ebres simples centrales 2-primaires ou 3-primaires 
co\"{\i}ncident.}

\medskip
Dans \cite{[CTGP]}, qui repose sur les travaux ant\'erieurs
de nombreux auteurs,
 on a \'etudi\'e de fa\c con syst\'ematique les propri\'et\'es
des corps de caract\'eristique z\'ero satisfaisant les hypoth\`eses  (1)
et (2), ou (1) et (2').
Les \'enonc\'es g\'en\'eraux suivants, extraits de \cite{[CTGP]},
s'appliquent donc aux corps de fonctions de deux variables
 sur un corps
alg\'ebriquement clos de caract\'eristique z\'ero.

\medskip

L'\'enonc\'e suivant (cf. \cite{[CTGP]}, Thm. 1.2)
 est un cas particulier de la
conjecture II de Serre (\cite{[S2]}, \cite{[S3]}).  Il rassemble des travaux
successifs de Merkur'ev-Suslin, Suslin, Bayer-Fluckiger et Parimala,
P.~Gille, Chernousov.

\begin{theo}[]\label{}
 Soit $K$ un corps   de caract\'eristique z\'ero
satisfaisant $(1)$  et $(2')$. 
 Soit $G$  un $K$-groupe alg\'ebrique  semi-simple simplement connexe.
 Si $G$ contient un facteur de type $E_8$, supposons en outre $\cd(K^{\ab})\leq 1$. Alors
 $H^1(K,G)=0$.
\end{theo}
On note ici $K^{\ab}$ l'extension ab\'elienne maximale de $K$.
La question de savoir si l'on a $\cd(K^{\ab})\leq 1$ pour un corps
$K$ de fonctions
de deux variables sur les complexes est ouverte.

\bigskip

Dire que sur un corps $K$ (de caract\'eristique z\'ero)  exposant et
indice co\"{\i}ncident pour les $K$-alg\`ebres simples
centrales
\'equivaut \`a dire que pour tout entier $n \geq 2$ l'application
bord $$H^1(K,PGL_n) \to H^2(K,\mu_n)$$
d\'eduite de la suite exacte
$$ 1 \to \mu_n \to SL_n \to PGL_n \to 1$$
est surjective.
On peut donc se poser des questions analogues avec d'autres
isog\'enies. Une longue analyse par type de groupe (le cas ${}^2A_n$ \'etant
particuli\`erement d\'elicat) m\`ene au r\'esultat suivant, qui requiert toute
l'hypoth\`ese (2),  et vaut donc pour les corps de fonctions de deux variables
sur les complexes.

\begin{theo}[\cite{[CTGP]}, Thm. 2.1]\label{}  Soit $K$  un corps
de caract\'eristique z\'ero satisfaisant  $(1)$ et $(2)$. Soit  $G$ un
$K$-groupe semi-simple  simplement
connexe, de centre $\mu$, de groupe adjoint $G^{\ad}$.
Si $G$ contient un facteur de type $E_8$, supposons en outre $\cd(K^{\ab})\leq 1$.

{\rm (i)}  L'application bord $H^1(K,G^{\ad}) \to H^2(K,\mu)$
associ\'ee
\`a la suite exacte
$$1 \to \mu \to G \to G^{\ad} \to 1$$
 est une bijection.

{\rm (ii)}  Si le groupe $G$ n'est pas purement de type $A$, alors il est isotrope.
\end{theo}

(Ce th\'eor\`eme g\'en\'eralise des r\'esultats classiques de Kneser sur
les corps locaux et sur les corps globaux.)

 \begin{theo}[cf. \cite{[CTGP]}, Thm. 4.5 ]\label{trivCTGP}  Soit
$K$ un corps de caract\'eristique z\'ero satisfaisant  $(1)$ et $(2')$. 
Soit
$G$ un
$K$-groupe semi-simple
simplement connexe. Si $G$ ne contient pas de facteur de type $E_8$,
ou si
l'on a en outre $\cd(K^{\ab})\leq 1$, alors $G(K)/R=1$.
\end{theo}

(Pour la notion de R-\'equivalence, voir  \cite{[G1]}.)

La d\'emonstration se fait en discutant chaque type simple. Il y a une
diff\'erence marqu\'ee entre le cas des groupes de type $A_n$ et les autres types simples. Pour les
autres types simples, la $K$-vari\'et\'e sous-jacente au $K$-groupe $G$
est $K$-birationnelle \`a un espace affine (d\'efini sur $K$),
ce qui n'est pas le cas en g\'en\'eral
pour les groupes de type $A_n$, comme le montre un exemple
de Merkur'ev.

 \bigskip

Le th\'eor\`eme \ref{trivCTGP} est utilis\'e dans la
d\'emonstration du th\'eor\`eme suivant, analogue d'un
 r\'esultat de P. Gille
(\cite{[G1]}, \cite{[G2]}).
\begin{theo}[\cite{[CTGP]}, Thm. 4.12 ]\label{} Soient $K$ un corps de
fonctions de deux variables sur un corps alg\'ebriquement clos de
caract\'eristique z\'ero et $G$ un $K$-groupe lin\'eaire connexe
sans facteur de type $E_8$. Le groupe $G(K)/R$ est un groupe
ab\'elien fini.
\end{theo}

(Ici encore la condition sur
$E_8$ pourrait \^etre omise si l'on savait \'etablir $\cd(K^{\ab})\leq 1$.)

 La technique des r\'esolutions
flasques de groupes lin\'eaires connexes permet  de donner une formule pour le
 groupe ab\'elien fini $G(K)/R$. Je renvoie pour cela \`a  \cite{[CTGP]},
\S 4.
 
\section{Quelques r\'esultats sur les corps de fonctions d'une variable
sur un corps $p$-adique}

 \begin{theo}[Saltman \cite{[Sa1]}]\label{Salt1} Soient  $k$ un corps
p-adique et
$K$ un corps de fonctions
  d'une variable
  sur $k$. Soit $l$ premier, $l \neq p$. Supposons $\mu_l \subset k$.
 Etant donn\'e un ensemble fini d'alg\`ebres simples centrales $A_i, i \in I$,
 chacune d'exposant $l$ dans le groupe de Brauer de $K$, il existe $f$ et $g$
 dans $K^*$ tels que l'extension $K(f^{1/l},g^{1/l})$ d\'eploie chacune des
alg\`ebres $A_i, i\in I$.
\end{theo}

L'outil de base de la d\'emonstration est le th\'eor\`eme suivant.

\begin{theo}[Tate, Lichtenbaum, Grothendieck]\label{H2} 
 Soit $O$ l'anneau des entiers d'un corps $p$-adique
et soit $Y/\Spec(O)$ un sch\'ema connexe, r\'egulier, propre et plat sur
$\Spec (O)$, de dimension relative $1$. Alors le groupe de Brauer de $Y$
est trivial.
\end{theo}

En utilisant la proposition \ref{1.3.1} on voit alors que 
 le {\og groupe de Brauer non ramifi\'e \fg} du corps
des fonctions $K$ d'une telle surface arithm\'etique $Y$
est trivial. En particulier si un \'el\'ement $A \in \Br(K)$
a la propri\'et\'e que pour tout anneau de valuation discr\`ete $R$
de corps des fractions $K$, de corps r\'esiduel $\kappa$,
 le r\'esidu $\delta(A) \in \kappa^*/\kappa^{*l}$ est trivial,
 alors $A=0 \in \Br(K)$. On notera que l'image
 du point ferm\'e $s_{R}$ de $\Spec (R)$ dans $Y$ peut \^etre soit
 le point g\'en\'erique d'une courbe soit un point ferm\'e de $Y$.

\medskip

La d\'emonstration du th\'eor\`eme \ref{Salt1} est  analogue 
dans son principe \`a celle
de \cite{[Ar3]}, \cite{[FS]}, \cite{[CTOP]} (o\`u l\`a on n'extrait qu'une
seule racine
$l$-i\`eme). On commence par \'ecrire le corps $K$ comme le corps
des fonctions d'un sch\'ema $X$ connexe, r\'egulier, propre et plat sur
$\Spec (O)$, de dimension relative $1$, choisi de telle sorte que
la r\'eunion des supports des diviseurs de ramification des
$A_i$ soit \`a croisements normaux stricts sur $X$; 
par \'eclatements convenables,
on peut assurer cela sur un sch\'ema r\'egulier excellent de dimension $2$
(Lipman \cite{[Li]}). 
On d\'etermine ensuite  des  fonctions $f,g \in K^*$
telles que les classes $A_i$ deviennent non ramifi\'ees sur le corps
$L=K(f^{1/l},g^{1/l})$, qui est le corps des fonctions d'un sch\'ema  $Y$ connexe, 
r\'egulier, propre et plat sur
$\Spec (O)$, de dimension relative $1$, muni d'un morphisme g\'en\'eriquement fini
$Y \to X$.
Le th\'eor\`eme \ref{H2} assure alors la nullit\'e
de ces classes dans $\Br(L)$. Pour choisir les fonctions $f$ et $g$,
on discute les images possibles dans $X$ des morphismes compos\'es
$s_{R } \to \Spec(R) \to Y \to X$, pour $R$ comme ci-dessus.
 Pour le choix (d\'elicat) des fonctions $f$ et $g$, je renvoie \`a \cite{[Sa1]}  et  \`a
\cite{[HVG]}.

\medskip
 Un argument alg\'ebrique simple
 permet de d\'eduire du th\'eor\`eme \ref{Salt1}
 l'\'enonc\'e suivant :

 \begin{coro}[Saltman \cite{[Sa1]}]\label{} Etant donn\'es un corps $K$
comme ci-dessus et
une alg\`ebre
simple centrale $A$ d'exposant $n$ premier \`a $p$, l'indice de $A$ divise 
$n^2$.
\end{coro}

  Des exemples dus \`a Jacob et Tignol  (\cite{[Sa1]}, voir aussi \cite{[KRTY]})
montrent que cette borne
  est en g\'en\'eral la meilleure. Sur le corps $K={\bf Q}_p(x)$,
  avec $p\neq 2$, si $a \in \Z_{p}$ est une unit\'e qui n'est pas 
  un carr\'e, alors 
   le produit tensoriel d'alg\`ebres de quaternions
  $(x,a)\otimes_{K} (x+1,p)$ est une alg\`ebre \`a division, donc d'indice 4
  mais d'exposant 2 (comparer avec l'exemple donn\'e apr\`es la proposition
  \ref{tignol}).

  Saltman vient r\'ecemment d'analyser ce type d'exemple.  Cela
lui a permis
  de montrer :
 \begin{theo}[Saltman \cite{[Sa2]}]\label{} Soient $k$ un corps p-adique et
$K$ un corps de fonctions
  d'une variable
  sur $k$. Soit $D$ une alg\`ebre \`a division sur $K$ d'indice $l$
premier, $l
\neq p$.
  Une telle alg\`ebre est cyclique~: $D$ contient un sous-corps commutatif
maximal
  cyclique sur $K$.
  \end{theo}

La d\'emonstration  passe par une  caract\'erisation des alg\`ebres de degr\'e
$l$ en termes de leur ramification, sur un mod\`ele r\'egulier propre convenable
de $K$ sur l'anneau des entiers du corps $p$-adique $k$. Elle utilise le
th\'eor\`eme \ref{H2}.

  \bigskip

 Un corps de fonctions d'une variable sur un corps local
d'\'egale caract\'eristique ${\bf F}_l((t))$ est un corps $C_3$. 
En particulier,
toute forme quadratique en au moins 9 variables sur un tel corps poss\`ede un
z\'ero non trivial. Il est donc naturel de poser la question :
si $K$ est un corps de fonctions d'une variable sur un corps
$p$-adique, toute forme quadratique en au moins 9 variables
poss\`ede-t-elle un z\'ero non trivial ?
  
 Merkur'ev remarqua le premier qu'en combinant
le th\'eor\`eme \ref{Salt1} (pour $l=2$) sur les corps de fonctions d'une variable
sur un corps $p$-adique (avec $p \neq 2$)   avec les r\'esultats g\'en\'eraux 
sur les formes quadratiques
obtenus gr\^ace \`a la K-th\'eorie alg\'ebrique (en fait uniquement
le th\'eor\`eme de Merkur'ev \cite{[M1]}), on obtient au moins une
borne sup\'erieure pour la dimension d'une forme quadratique
anisotrope d\'efinie sur un tel corps. Son r\'esultat fut am\'elior\'e par
Hoffmann et Van Geel  \cite{[HVG]} puis par Parimala et Suresh  :

 \begin{theo}[Parimala/Suresh \cite{[PS1]}]\label{}  Soit $K$ un corps de
fonctions d'une variable sur
  un corps $p$-adique, avec $p \neq 2$. Toute forme quadratique 
en au moins $11$ variables 
sur $K$
  poss\`ede un z\'ero non trivial.
  \end{theo}

Avant d'indiquer le principe de la d\'emonstration,
commen\c cons par quelques rappels. Le groupe de Witt $WK$
d'un corps $K$ de caract\'eristique diff\'erente de $2$  est par d\'efinition l'ensemble des
classes d'isomorphie de formes quadratiques non d\'eg\'en\'er\'ees, muni de la
somme directe orthogonale et quotient\'e par la classe du plan
hyperbolique standard
$x^2-y^2$. Le produit tensoriel des formes quadratiques lui donne
une structure d'anneau. L'id\'eal des classes de formes de rang pair
est not\'e $IK \subset WK$. L'id\'eal $I^nK$ est engendr\'e additivement
par les $n$-formes de Pfister. L'intersection de tous
les $I^nK$ est r\'eduit \`a z\'ero. Plus pr\'ecis\'ement, toute forme
anisotrope appartenant \`a $I^nK$ a une dimension au moins \'egale
\`a $2^n$ (Arason-Pfister). La conjecture de Milnor \'etablie
par Voevodsky implique (Orlov-Vishik-Voevodsky)
l'existence d'isomorphismes
$I^nK/I^{n+1}K \simeq H^n(K,\Z/2)$. Dans ces isomorphismes,
une $n$-forme de Pfister correspond \`a un $n$-symbole.

On peut alors facilement \'etablir l'\'enonc\'e g\'en\'eral suivant.

 \begin{prop}[\cite{[Kh]}, \cite{[PS2]}]\label{genKPS} Soit $K$ un corps de
caract\'eristique diff\'erente de $2$.
Soit $N>0$ un entier tel que $H^{N+1}(K,\Z/2)=0$; on suppose que pour tout
entier $n$ avec $1 \leq n \leq N$ il existe un entier $\lambda_n(K)$ tel
que tout \'el\'ement de $H^n(K,\Z/2)$ est somme d'au plus $\lambda_n(K)$
symboles. Alors il existe un entier positif $u(K)$ tel que toute
forme quadratique en strictement plus de $u(K)$ variables ait un z\'ero 
non trivial.
\end{prop}
\noindent{\sc Preuve} --- L'hypoth\`ese et les rappels ci-dessus 
impliquent d'une part que  $I^{N+1}K=0$
(et  $H^m(K,\Z/2)=0$ pour tout $m>N$), d'autre part
que tout \'el\'ement de $I^nK$ peut s'\'ecrire comme une somme
orthogonale de $\lambda_n(K)$ $n$-formes de Pfister
et d'une forme appartenant \`a $I^{n+1}K$. Ainsi dans $WK$ toute
forme quadratique est repr\'esent\'ee par une
 forme de rang au plus $1+\sum_{n=1}^N 2^n \lambda_n(K)$.
Ceci implique que toute forme quadratique de rang strictement plus
grand que $1+\sum_{n=1}^N 2^n \lambda_n(K)$ est isotrope.
\cqfd

La majoration obtenue peut \^etre am\'elior\'ee  (\cite{[Kh]} Prop. 1.2.d) :
$$u(K) \leq 2+\sum_{n=2}^N (2^n-2) \lambda_n(K).$$

Dans le cas qui nous int\'eresse ici, le corps $K$
est de dimension cohomologique $3$ par des arguments g\'en\'eraux,
en particulier  $H^4(K,\Z/2)=0$.
On a la borne \'evidente $\lambda_1(K)=1$ (valable sur tout corps).
 Le th\'eor\`eme \ref{Salt1} pour $l=2$  et un r\'esultat  d'Albert 
 (\cite{[A]},  Chap. XI, Thm. 2; \cite{[Ar2]}, Thm. 5.5)
 impliquent
que tout \'el\'ement de $H^2(K,\Z/2)$ est la classe d'un produit
tensoriel de deux alg\`ebres de quaternions, i.e. est la somme
  de deux symboles.
On a donc 
$\lambda_2(K)=2$. Il reste \`a majorer $\lambda_3(K)$.
Partant du r\'esultat de Saltman, par des manipulations
alg\'ebriques, Hoffmann et Van Geel \cite{[HVG]} \'etablissent
$\lambda_3(K)
\leq 4$, puis $u(K) \leq 22$.
 Par un argument de g\'eom\'etrie arithm\'etique,
Parimala et Suresh montrent $\lambda_3(K)=1$~:
  \begin{theo}[\cite{[PS1]}]\label{H3}
  Pour $K$ comme ci-dessus, toute classe dans $H^3(K,\Z/2)$
est repr\'esentable par un seul symbole $(a)\cup (b) \cup (c)$.
\end{theo}
 La majoration  g\'en\'erale ci-dessus donne alors imm\'ediatement $u(K) \leq 12$.
Un travail arithm\'etique plus fin permet \`a Parimala et Suresh d'obtenir la  borne $u(K) \leq 10$.
 L'outil fondamental pour la d\'emonstration de \ref{H3} est 
 le r\'esultat suivant (analogue pour le groupe $H^3$
  du r\'esultat \ref{H2}), qui est un cas
   particulier d'un th\'eor\`eme de  K. Kato~:
 \begin{theo}[\cite{[Kt]}]\label{}
  Soient $k$ un corps $p$-adique avec $p\neq 2$ et  $K$
  un corps de fonctions d'une variable sur $k$. Le groupe de cohomologie
  non ramifi\'e $H^3_{\nr}(K,\Z/2)$ est nul.
  \end{theo}

  Parimala et Suresh partent d'un \'el\'ement quelconque  $\alpha \in
H^3(K,\Z/2)$.
 En consid\'erant un mod\`ele r\'egulier propre de $K$ au-dessus de 
l'anneau
 des entiers de $k$, ils montrent comment trouver un \'el\'ement  $f \in
K^*$
 tel que $\alpha$ devienne non ramifi\'e dans $H^3(K(\sqrt{f}),\Z/2)$, 
donc nul
 dans $H^3(K(\sqrt{f}),\Z/2)$
 par le th\'eor\`eme de Kato appliqu\'e \`a $K(\sqrt{f})$. Le th\'eor\`eme 
de Saltman
 assure alors que $\alpha$ est une somme d'au plus 2 symboles. Un travail
 arithm\'etique plus pr\'ecis montre que $\alpha$ est repr\'esent\'e par un seul symbole.
 Le th\'eor\`eme \ref{H3}  r\'epond \`a une question de Serre \cite{[S3]} :
 pour $K$ comme ci-dessus et $G$ le $K$-groupe simple d\'eploy\'e de type $G_{2}$,
 la fl\`eche naturelle $H^1(K,G)  \to H^3(K,\Z/2)$ est une bijection.

\end{document}